\documentclass{article}

\usepackage{preamble}

\title{Pointed quandles of linkoids}
\author{
Neslihan Gügümcü\footnote{\href{mailto:neslihangugumcu@iyte.edu.tr}{neslihangugumcu@iyte.edu.tr}, Department of Mathematics, Izmir Institute of Technology, Urla Izmir 35430, Turkey}
\and 
Runa Pflume \footnote{\href{mailto:runa@pflume.de}{runa@pflume.de}, Institut of Mathematics, Georg-August University Göttingen, Bunsenstraße 3-5
37073 Göttingen, Germany }}

\begin{document}
	\renewcommand{\dobib}{} 
	\maketitle

\begin{abstract}
In this paper we define the fundamental quandle of knotoids and linkoids and prove that it is invariant under the under forbidden-move and hence encodes only the information of the underclosure of the knotoid. \\
We then introduce $n$-pointed quandles, which generalize quandles by specifying $n$ elements as ordered basepoints. This leads to the notion of fundamental pointed quandles of linkoids, which enhances the fundamental quandle. 
Using $2$-pointed quandle allows us to distinguish 1-linkoids with equivalent under-closures and leads to a couple of 1-linkoid invariants. \\
We also generalize the notion of homogeneity of quandles to $n$-homogeneity of quandles. We classify all $\infty$-homogeneous, finite quandles. 
\end{abstract}

\section{Introduction}\label{introsection}

Knotoids can be considered as knot diagrams in a plane or in the 2-sphere $S^2$ with two distinct open endpoints which cannot be moved over or under other strands in the diagram. Knotoids were introduced by Turaev in 2010 [Tur10], and since then, there have been numerous studies investigating properties and invariants of knotoids, see for example [GK17], [Mol22], [BBHL19], [MNK]. Turaev also introduced 1-linkoids (under the name multi-knotoids) as knotoids with additional closed (knot) components. Generalizations of 1-linkoids that contain a number of knotoid and knot components were studied in [GG22]. 

Quandles, on the other hand, are algebraic structures admitting a number of properties that give rise to a coloring of an oriented knot or link diagram, which is preserved under Reidemeister moves of knots, so that the coloring induced is indeed an invariant of oriented knots or links. Quandles were introduced independently by Joyce in [Joy82] and Matveev in [Mat84] (under the name \textit{distributive groupoids}). Both proved that knot quandles can distinguish all oriented knots up to mirror image with reversed orientation. This implies that quandles are as complex as knots themselves. The reader may find more about quandles in [EMRL10], [Nel02] and [BLRY10] and more on  methods to distinguish different knots and links using quandles in [CJK+01], [CN18] and [CESY14].

Quandle colorings of 1-tangles were studied in [CSV16] and [CDS16, Chapter 3]. 1-tangles can be understood as knot-type knotoids that are a specific type of knotoids with both endpoints lying in the same region of the diagram. Generalizations of quandle colorings, such as biquandle colorings and shadow quandle colorings were also studied for knotoids in [GN18] and [Caz22], respectively.  There is, however, still only very little research on  knotoids and linkoids in relation with quandles.

In this paper we define the fundamental quandle of a linkoid and study its
basic properties. We show that the fundamental quandle is an invariant of linkoids but it remains unchanged also under moving an endpoint under another arc, which is in fact a forbidden move for linkoids.  This move is forbidden since it may transform the linkoid to a non-equivalent one. To enhance the fundamental quandle for linkoids, we introduce a generalization of quandles, called
\textit{$n$-pointed quandles}. These are quandles with n ordered basepoints. This generalization leads to the new concept of fundamental pointed quandle of a linkoid. We provide examples of knotoids that cannot be distinguished by using quandles but can be distinguished by using pointed quandles. We also study homogeneity of quandles. A quandle $X$ is homogeneous if the group of quandle automorphisms $Aut(X)$ of the quandle acts transitively on it. With our new definition of n-pointed quandles, we introduce the notion of quandles being $n$-homogeneous. We classify all finite quandles that are $k$-homogeneous, where $k$ is the cardinality of the quandle. This paper is an excerpt of the second author's thesis \cite{Pfl23}, supervised by the first author.

Let us now introduce the organization of the paper. In Section \ref{preliminariessection} we recall preliminary notions and fundamental theorems of the theory of knotoids and linkoids, and the theory of quandles. In Section \ref{fundamentalsection}, we introduce and study the fundamental quandle of linkoids. We show that the fundamental quandle is an invariant of linkoids. We observe that the fundamental quandle cannot distinguish two linkoids which are related to each other by under-forbidden. In particular, we show that if two knotoids represent the same knot then their fundamental quandles are isomorphic. In Section \ref{pointedsection} we  add extra structure to quandles, and introduce $n$-pointed quandles. $n$-pointed quandles yield to the fundamental pointed quandle of a linkoid. We study the fundamental pointed quandles of linkoids and show that the fundamental pointed quandle is a stronger invariant than the fundamental quandle for 1-linkoids. We also classify all $\infty$-homogeneous, finite quandles and give a lower bound for the number of pointed quandle with a given underlying finite quandle, up to isomorphism. Finally, we introduce the pointed quandle counting invariant of a linkoid and study the quandle counting matrix of a 1-linkoid.  In Section \ref{discussionsection} we note down further problems.

\section{A review of linkoids and quandles}\label{preliminariessection}

\subsection{Knotoids and Linkoids}
    In this section we recall fundamental notions from the theory of knotoids and linkoids. 

    \begin{defn}    
    An oriented \emph{n-linkoid diagram} in $S^2$ is a generic immersion of $n \geq 0$ unit intervals $[0,1]$ oriented from $0$ to $1$ and a number of oriented unit circles $S^1$ into $S^2$ with finitely many transverse double points. Each such double point is endowed with over/under-crossing data. We call a component \textit{open}, if it is the image of $[0,1]$ and \textit{closed} if it is the image of $S^1$. The image of $0$ in an open component is named as the \emph{leg} (or tail) and the image of $1$ is named as the \emph{head} of the component. 
    
    In particular, we call a 1-linkoid diagram with no closed components a \textit{knotoid diagram}, and a 0-linkoid with only one closed component a \textit{knot} and more closed components a \textit{link}. A \textit{full linkoid} is an n-linkoid with no closed components, where $n > 1$.
    \end{defn}

\begin{figure}[ht]
	\begin{center}
        \begin{subfigure}[b]{0.45\textwidth}
            \centering
            \includegraphics[height=3cm]{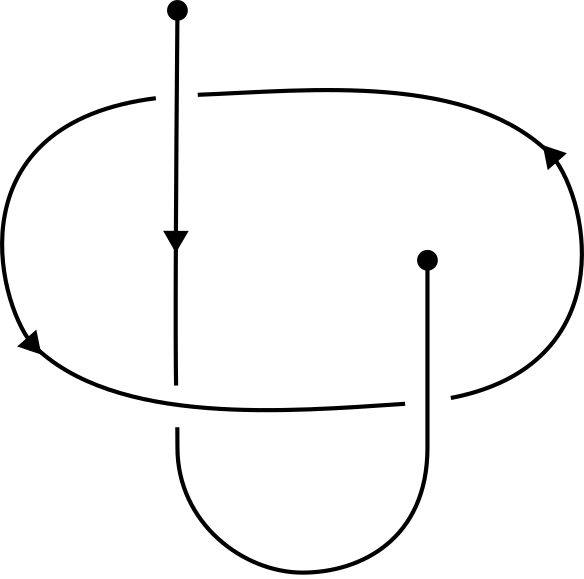}
            \caption{A 1-linkoid diagram with two components}
        \end{subfigure}
        \begin{subfigure}[b]{0.45\textwidth}
            \centering
            \includegraphics[height=3cm]{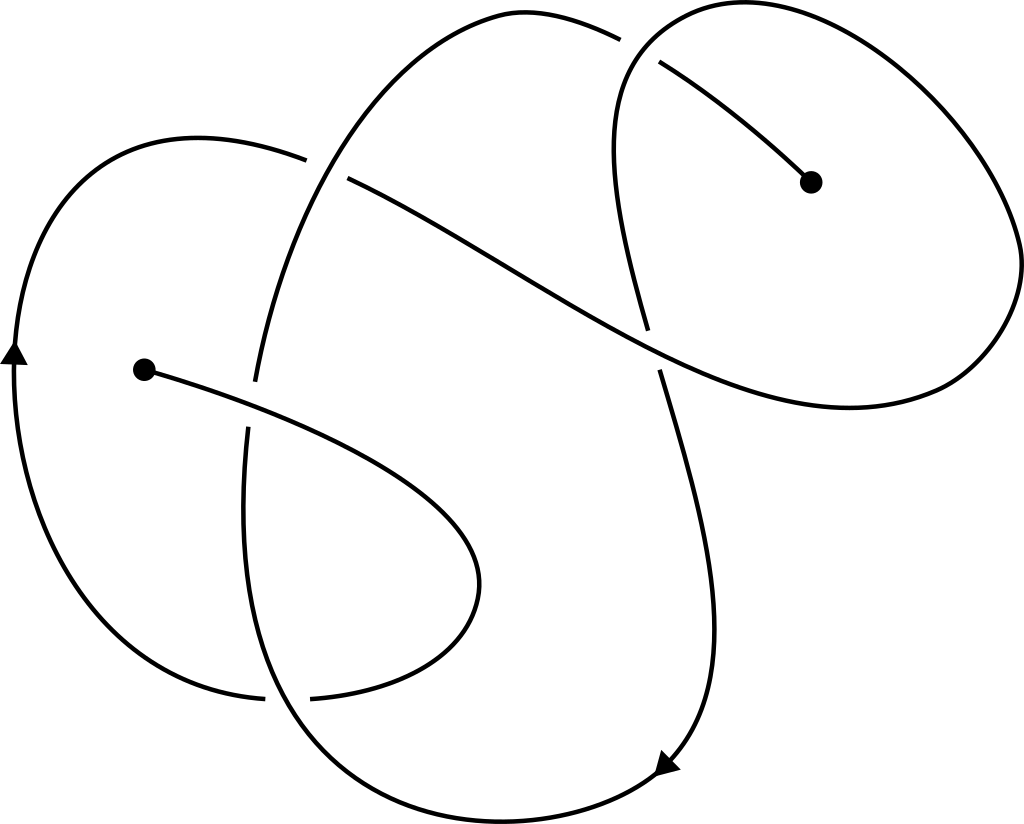}
            \caption{A knotoid diagram}
            \label{528}
        \end{subfigure}
    \end{center}
\end{figure}

 There are two types of crossings with respect to the orientation, called positive and negative crossing as specified in Figure \ref{crossingtypes} below.

\begin{figure}[ht]
	\begin{center}
        \begin{subfigure}[b]{0.3\textwidth}
            \centering
            \includegraphics[width=0.5\textwidth]{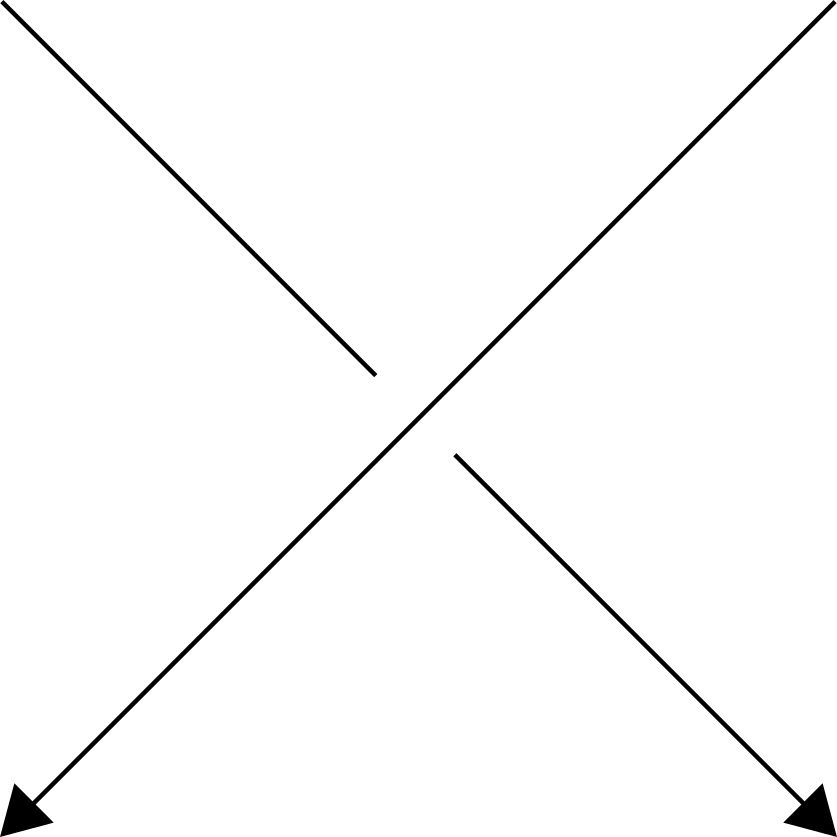}
            \caption{Positive crossing}
            \label{positive_crossing}
        \end{subfigure}
        \begin{subfigure}[b]{0.3\textwidth}
            \centering
            \includegraphics[width=0.5\textwidth]{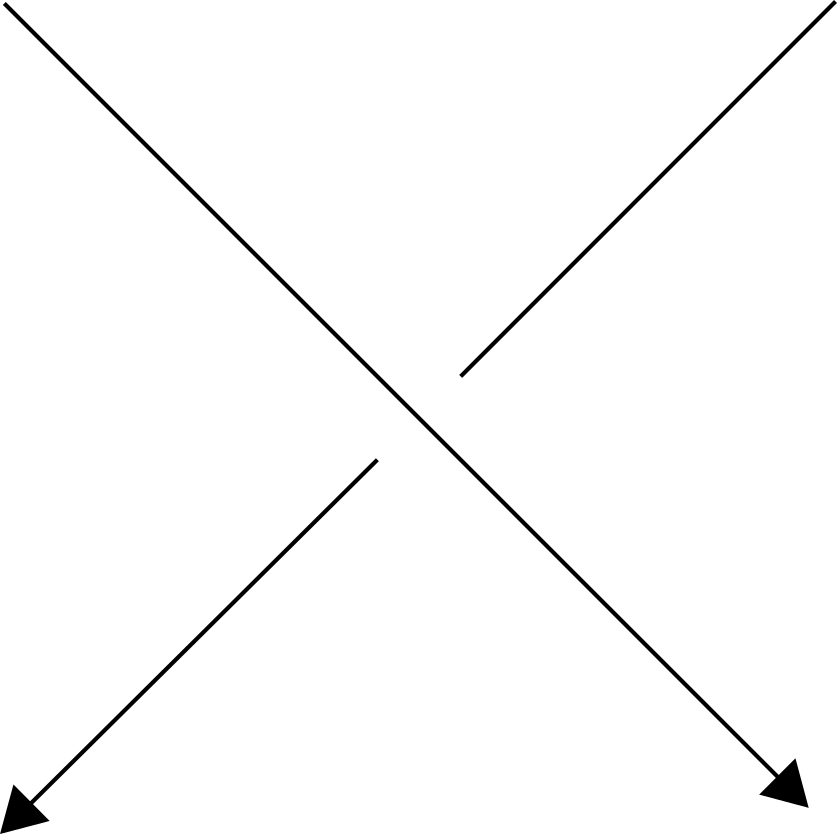}
            \caption{negative crossing}
            \label{negative_crossing}
        \end{subfigure}
        \caption{}
        \label{crossingtypes}
	\end{center}
\end{figure}

\begin{defn}

    An n-linkoid that is equivalent to a linkoid with zero crossings is called a \emph{trivial} n-linkoid. The trivial knot is called the \emph{unknot}.
\end{defn}

\begin{rmk}
\begin{enumerate}~\\
    \item 
In Turaev's paper \cite{Tur10}, a 1-linkoid is called a multi-knotoid. A full linkoid is called a linkoid in \cite{GG22}.
\item 

    Linkoids can in general be defined as generic immersions into any orientable surface. We will only consider spherical linkoids in this paper. More general surfaces lead to the notion of a virtual knot(oid) and link(oid) diagrams. These were first introduced in \cite{Kau98} and further studied for instance in \cite{GK17}, \cite{KR01} or \cite{FJK04}.
    \end{enumerate}
\end{rmk}

\begin{figure}[ht]
	\begin{center}
        \begin{subfigure}[b]{0.3\textwidth}
            \centering
            \includegraphics[height=2.5cm]{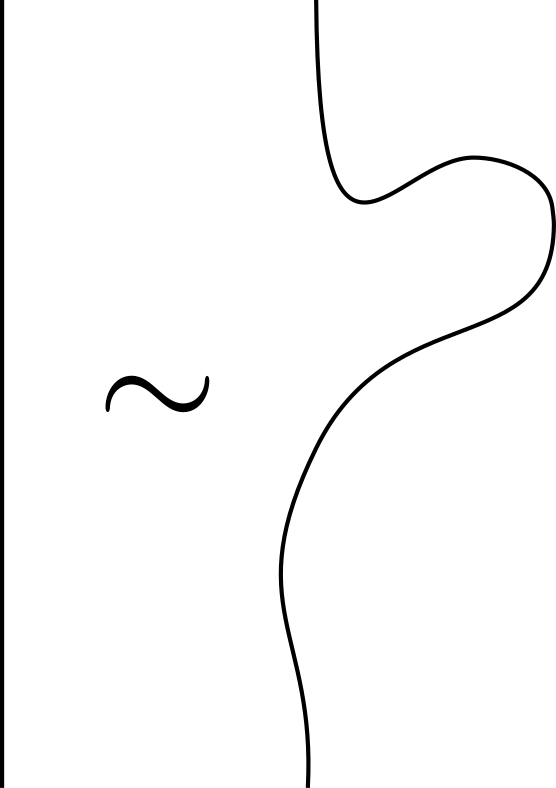}
            \caption{R0-move}
            \label{r0}
        \end{subfigure}
        \begin{subfigure}[b]{0.3\textwidth}
            \centering
            \includegraphics[height=2.5cm]{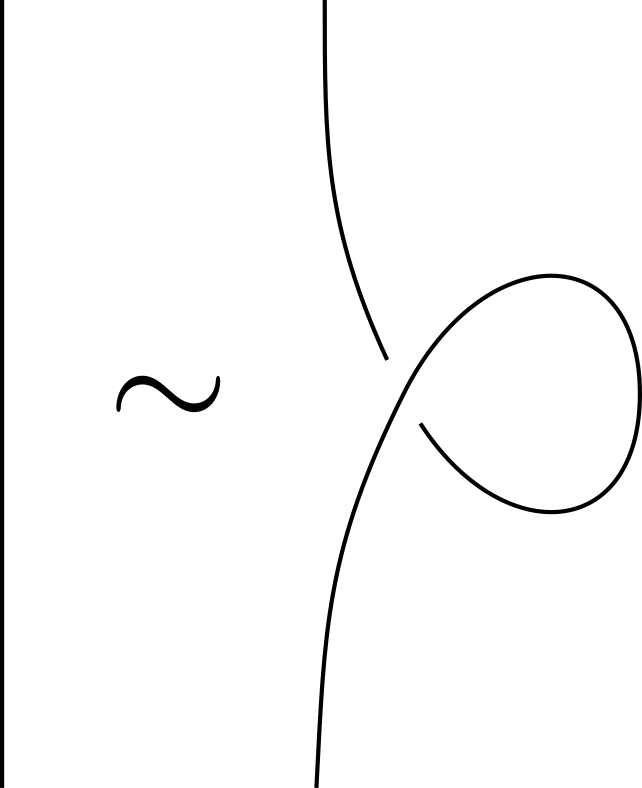}
            \caption{R1-move}
            \label{r1}
        \end{subfigure}
        \begin{subfigure}[b]{0.3\textwidth}
            \centering
            \includegraphics[height=2.5cm]{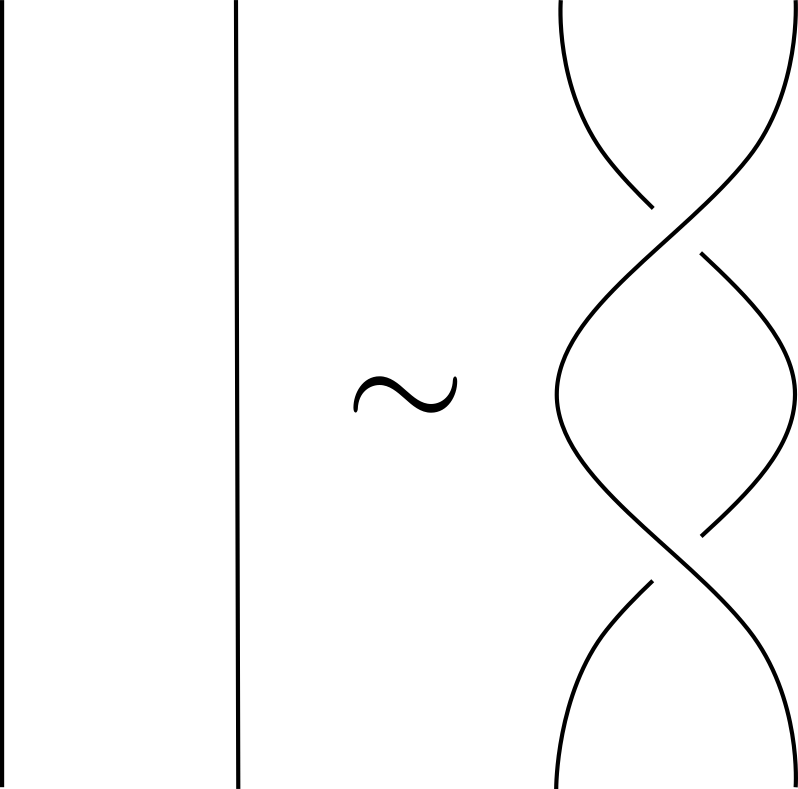}
            \caption{R2-move}
            \label{r2}
        \end{subfigure}
        \begin{subfigure}[b]{0.49\textwidth}
            \centering
            \includegraphics[height=2.5cm]{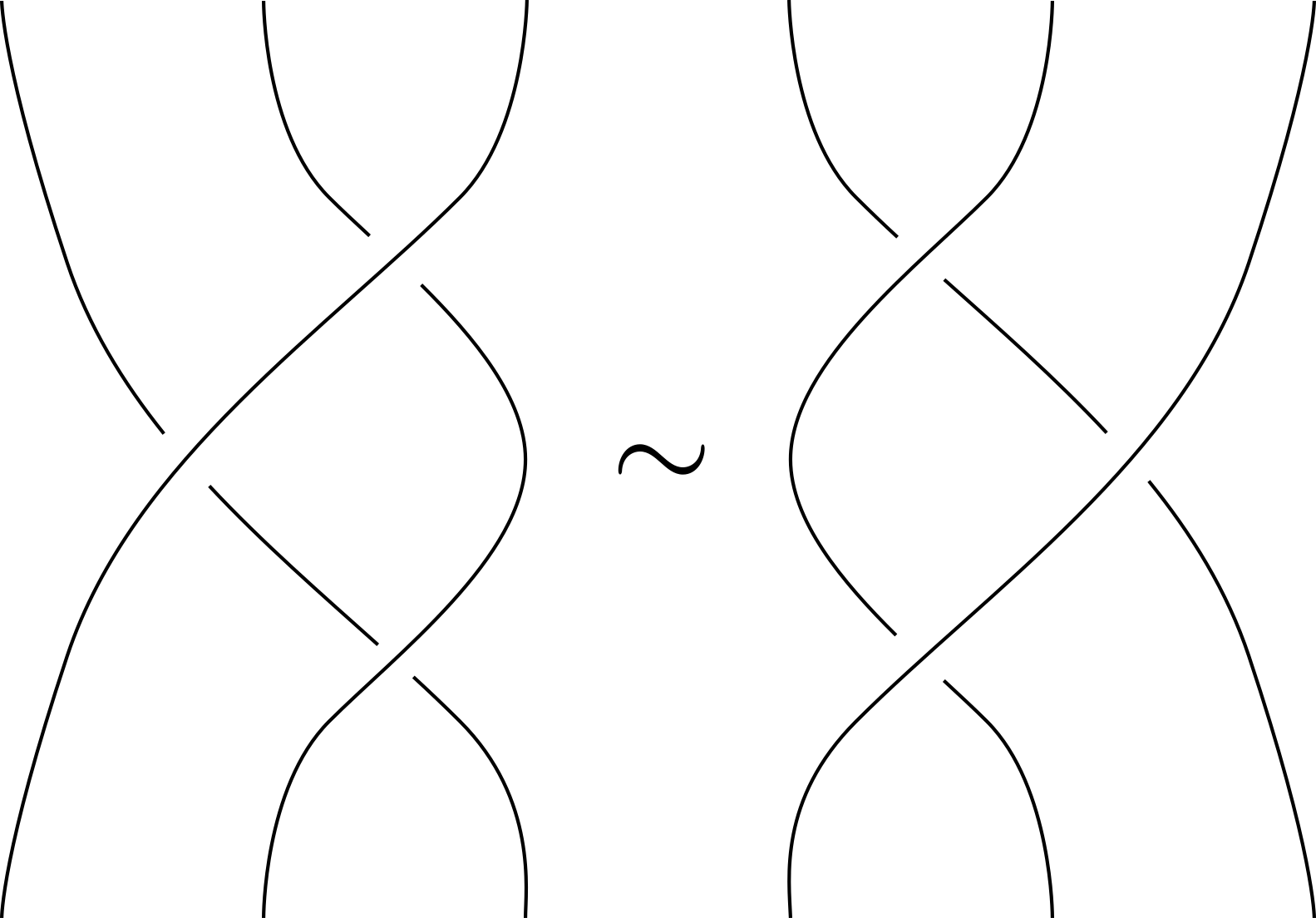}
            \caption{R3-move}
            \label{r3}
        \end{subfigure}
        \begin{subfigure}[b]{0.49\textwidth}
            \centering
            \includegraphics[height=2.5cm]{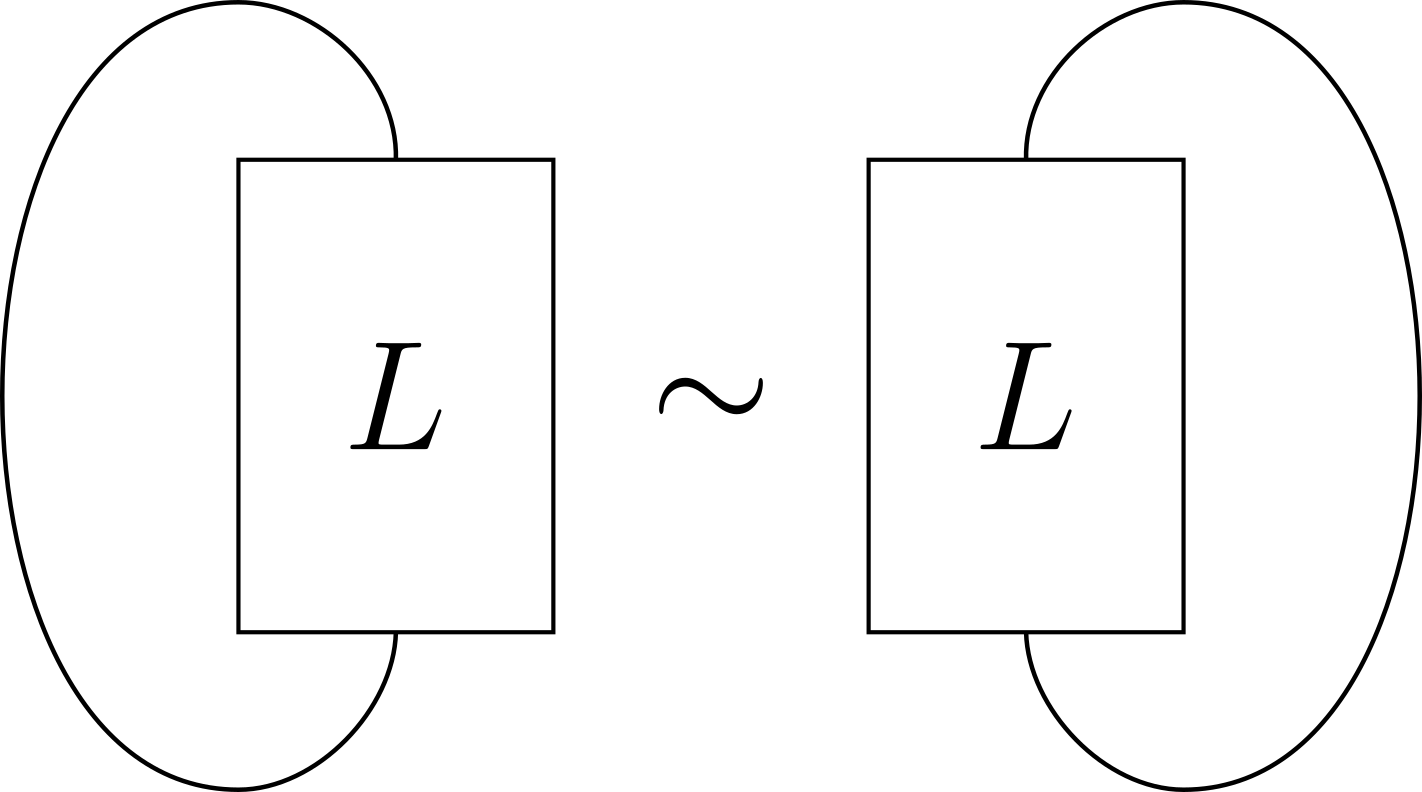}
            \caption{S-move}
            \label{smove}
        \end{subfigure}
        \caption{Reidemeister moves (a) - (d) and the spherical move (e)}
        \label{rmoves}
	\end{center}
\end{figure}

\begin{defn}
    Two linkoid diagrams $L, L'$ are \emph{equivalent}, if one can be moved into the other by a finite sequence of local oriented Reidemeister moves R0 - R3 depicted in Figure \ref{rmoves}. These moves happen away from the endpoints and cannot move endpoints over or under a strand. We write $L \sim L'$. The equivalence classes of these diagrams are called \emph{linkoids.}
\end{defn}

\begin{figure}[ht]
    \begin{center}
        \includegraphics[height=2cm]{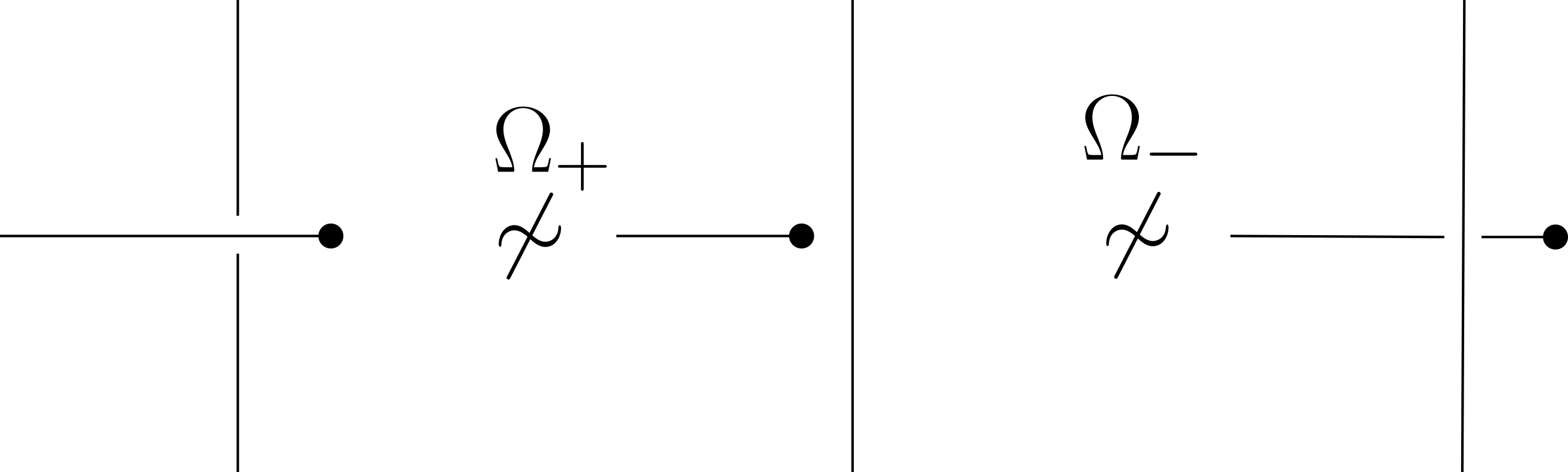}
        \caption{The over and under forbidden under}
        \label{forbidden}
    \end{center}
\end{figure}

\begin{rmk}
    In Figure \ref{rmoves} we depict Reidemeister moves that take place on an unoriented linkoid diagram. The reader can verify that there are in total four oriented versions of R1, four oriented versions of R2, and eight oriented versions of R3 moves. In \cite{pol10}, a minimal generating set is shown to contain four moves consisting of two oriented R1, one oriented R2 and one oriented R3 move. 
\end{rmk}
\begin{defn}
Pulling an endpoint over or under an adjacent strand may change the type of an n-linkoid. In fact, allowing an endpoint to be pulled both over and under other strands clearly turns any linkoid diagram to a trivial linkoid diagram. These moves, also depicted in Figure \ref{forbidden} are called the \emph{forbidden moves} and are denoted by $Ω_+$ and $Ω_-$, respectively. \\
    
\end{defn}

Given a 1-linkoid diagram $L$. We can obtain a link from $L$ by connecting the head to the leg of $L$ with an embedded arc that intersects $L$ only transversely at a finite set of points all of which are endowed with either under or over information. We call this arc a \emph{shortcut} of the 1-linkoid. If the shortcut goes under every strand of $L$ to connect its endpoints, the resulting link, denoted by $L_-$ is called the \emph{under closure} of $L$. Note that every choice of the shortcut gives an equivalent closure up to Reidemeister moves. If the shortcut goes over every strand of $L$, the resulting link, denoted by $L_+$ is called the \emph{over closure}. \\
We say that the 1-linkoid diagram $L$ \emph{represents} the link $L_-$.

Conversely, let $\tilde{L}$ be an oriented link diagram. Deleting an underpassing strand (a strand that goes under each crossing or containing no crossing at all) from a component of $\tilde{L}$ clearly results in a 1-linkoid diagram. It is also clear that the under closure of this resulting 1-linkoid diagram with the deleted strand gives us back $\tilde{L}$.

\begin{defn}

A 1-linkoid diagram is called a \textit{link-type} if the endpoints of the open component lie in the same region of the diagram. Specifically, if we have a knotoid diagram with two endpoints lying in the same region, it is called a \textit{knot-type} knotoid.
    
\end{defn}
Note that if $L$ is a link-type linkoid then $L_- = L_+$ and we say $L$ is \textit{of type} $L_-$. \\

Lemma \ref{ltlrep} was given in \cite[Section 2.2]{Tur10} for knotoids. The statement can be generalized to 1-linkoids directly. 
\begin{lem}\label{ltlrep}
    Two 1-linkoid diagrams $L$ and $L'$ represent the same link, if they can be transformed into each other by a finite sequence of Reidemeister moves and under forbidden moves $Ω_-$.
\end{lem}

\begin{prop}\label{knotktkequiv}[Tur10]
    There is a one-to-one correspondence between knots and knot-type knotoids.
\end{prop}

Because we can transform any knotoid into a knot-type knotoid using under forbidden moves, Proposition \ref{knotktkequiv} together with Lemma \ref{ltlrep} proves the next corollary.

\begin{cor}\label{omegacheat}
   Two knotoids represent the same knot if and only if they can be transformed into each other by a finite sequence of Reidemeister moves and under forbidden moves $Ω_-$.
\end{cor}

Note that Corollary \ref{omegacheat} is not true for links and 1-linkoids because by deleting an underpassing strand from different components of a link may result in non-equivalent 1-linkoid diagrams whose under closures are the same link.

\subsection{Quandles}

\begin{defn}\label{quandle}
    A \emph{quandle} is a set $ (X, \q) $ equipped with an operation $ \q : X \times X \to X $ such that
    \begin{enumerate}[label=(\arabic*)]
        \item $ x \q x = x$ for all $ x \in X $. (idempotence)
        \item For any $y \in X$, the map $β_y : X \to X$ defined as $β_y (x) = x \q y$ is bijective. (right invertibility) 
        \item $ (x \q y) \q z = (x \q z) \q (y \q z) $ for all $ x,y,z \in X $ (right self-distributivity).
     \end{enumerate}
\end{defn}
We often only write $X$ for the quandle $(X,\q)$, if the structure is clear. We write $x \q^{-1} y \coloneq β_y^{-1}(x)$ and see immediately 
\[(x \q y) \q^{-1} y= x = (x \q^{-1} y) \q y \]
 for all $x,y \in X$. We think of $x \q y$ as $y$ “acting” on $x$. For $x \q y$ we read it as “$x$ quandle $y$” or “we quandle $x$ with $y$”. A quandle is in general neither commutative nor associative. This means it is important to write parentheses.

\begin{exmp}~
\begin{itemize}
    \item Let $X$ be any set and $x \q y = x$ for all $x \in X$. This defines a quandle structure. It is called the \emph{trivial quandle}. If $|X| = n$ is finite, it is denoted by $T_n$. In particular a trivial quandle can have any number of elements, unlike for example the trivial group, which has only one element.
    \item Let $n \in ℕ $ and $X = ℤ/nℤ $ or $X = ℤ$. We can define a quandle structure on $X$ by $x \q y \coloneq 2y - x ( \text{mod }n)$.
    This quandle is called the \emph{dihedral quandle} and denoted as $R_n$.
    \end{itemize}
\end{exmp}
 
 \begin{defn}
    A \emph{quandle homomorphism} between two quandles $(X, \q_X)$ and $(Y, \q_Y)$ is a map $f \colon X \to Y$ such that 
    \[
        f(x_1 \q_X x_2) = f(x_1) \q_Y f(x_2).
    \]
     We denote the category of quandles by $\Qnd$.
     A bijective quandle homomorphism is called an \emph{isomorphism}.
\end{defn}


\begin{defn}\label{autodef}
    For a quandle $X$, the group of all isomorphisms $X \to X$ is called the \emph{automorphism group of $X$} and denoted by $Aut(X)$. \\
    The subgroup generated by the isomorphisms $β_y :x \mapsto x \q y$ for every $y \in X$ is called the group of \emph{inner automorphisms of $X$} and denoted by $Inn(X)$. 
\end{defn}

\begin{defn}\label{componentdef}
	The \emph{algebraic components} of a quandle $Q$ are the orbits under action of the inner automorphism group $Inn(Q)$.
\end{defn}

\begin{defn}\label{connecteddef}
    A quandle $X$ is \emph{connected}, if it has only one algebraic component. So for each $x, y \in X$ there are $x_1, …, x_n \in X$ and $ε_1, …, ε_n \in \{ -1, 1 \} $ such that 
    \[
        (…((x \q ^{ε_1} x_1) \q ^{ε_2} x_2) … ) \q^{ε_n}  x_n = y.
    \]
    This is equivalent to say that a quandle $X$ is \emph{connected} if $Inn(X)$ acts transitively on $X$. So for each $x,y \in X$ there is an $f \in Inn(X)$ with $f(x) = y$. 
\end{defn}

\begin{defn}\label{homogeneousdef}
    A quandle $X$ is \emph{homogeneous} if $Aut(X)$ acts transitively on $X$. So for each $x,y \in X$ there is an $f \in Aut(X)$ with $f(x) = y$. 
\end{defn}

\begin{defn}\label{faithfuldef}
    A quandle $X$ is \emph{faithful} if the map $x \mapsto β_x$ is injective. 
\end{defn}

\begin{figure}[ht]
   \begin{center}
      \begin{subfigure}[b]{0.3\textwidth}
         \centering
         \includegraphics[height=2.5cm]{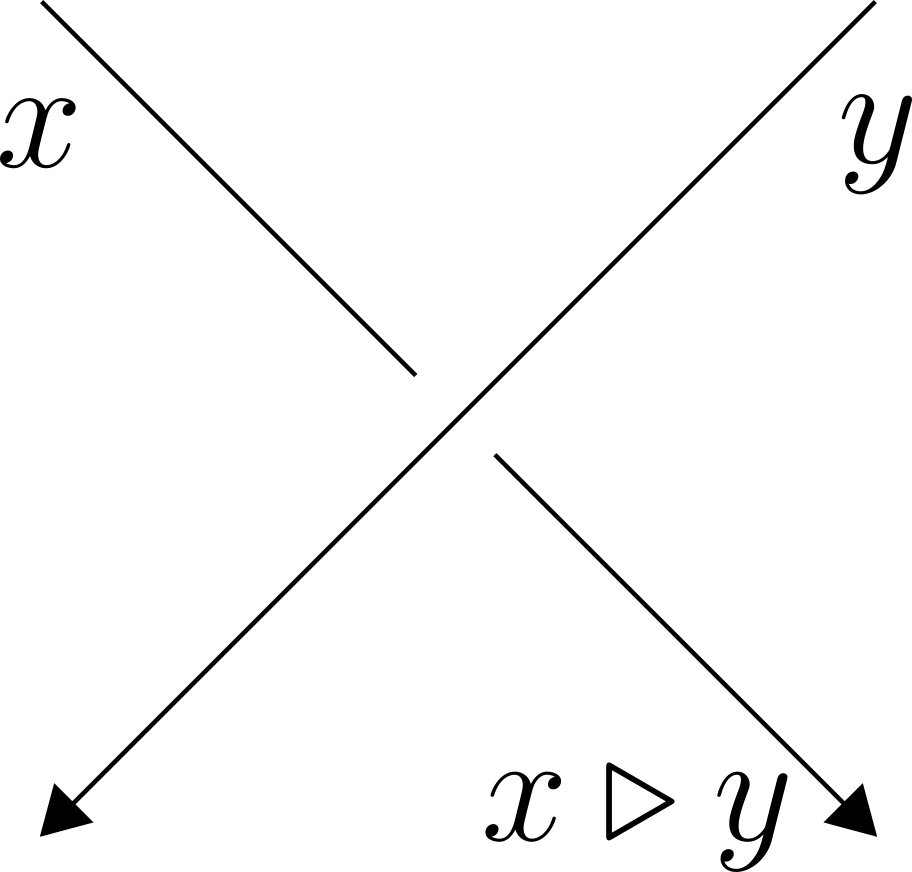}
      \end{subfigure}
      \begin{subfigure}[b]{0.3\textwidth}
         \centering
         \includegraphics[height=2.5cm]{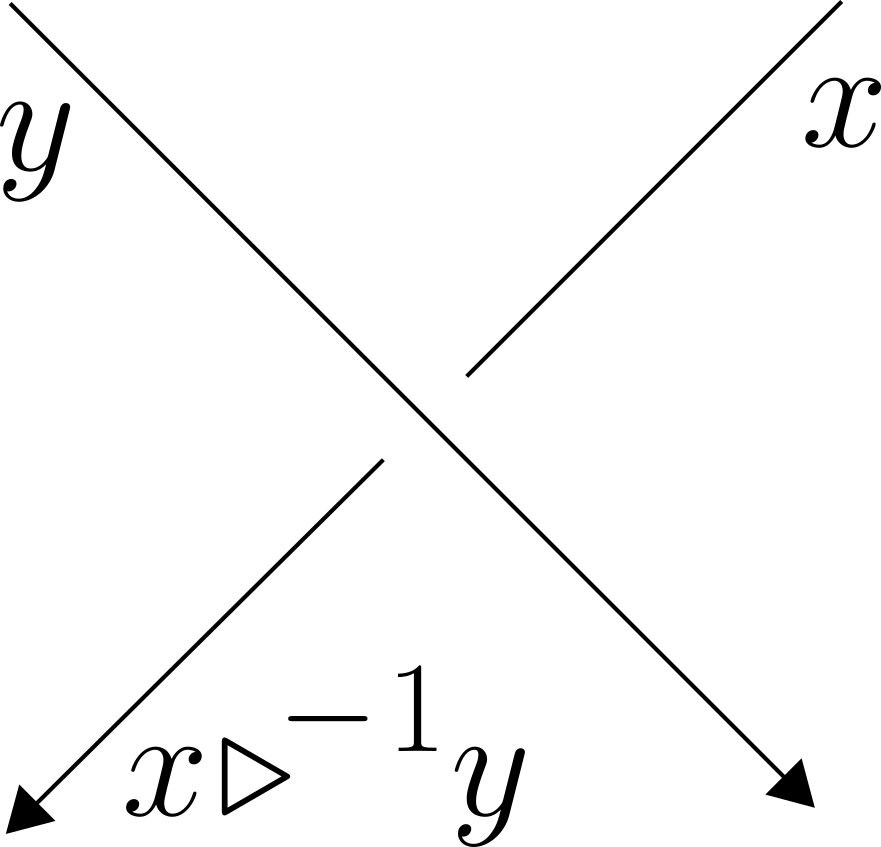}
         \label{quandlecond2}
      \end{subfigure}
      \caption{Quandle relation given on a crossing}
      \label{quandlecond}

   \end{center}

\end{figure}

Although, quandles are defined as algebraic structures with certain axioms, they can be considered as coloring sets for oriented knots and links. Precisely, a quandle coloring of an oriented knot or link diagram is the assignment of elements of a given quandle to arcs of the diagram so that the relations given in Figure \ref{quandlecond} are satisfied at positive and negative crossings of the diagram. With this point of view, quandle axioms are indeed the conditions induced by the quandle colored diagram with the requirement that the coloring remains invariant under oriented Reidemeister moves.  The reader is referred to \cite{en15} for the verification of this. We discuss below the fundamental quandle which is a powerful invariant of oriented knots and links, and use it to formally define quandle coloring. 


\begin{defn}[{\cite[Chp. 15]{joy82}, \cite{Mat84}}]\label{quandledeflink}
Let $L$ be an oriented link diagram and $A(L)$ the set of arcs in the diagram. Then the \emph{fundamental quandle} of $L$ is defined as 
\[
   Q(L) \coloneq Q\langle x \in A(L) \quad | \quad r_τ \text{ for all crossings } τ \rangle,
\]
where the quandle consists of words in $A(L)$ modulo the quandle axiom relations and the relations $r_τ$ given by each crossing $τ$ as in Figure \ref{quandlecond}. 
\end{defn}

\begin{thm}
    The fundamental quandle is invariant under Reidemeister moves and hence only depends on the link. 
\end{thm}

For a given oriented link $L$ we define the reversed link $r(L)$ as $L$ with reversed orientation of every component and the \emph{mirror link} $m(L)$ where we change every over crossing to an under crossing and vice versa. We write $rm(L)$ for the reversed mirror of $L$. \\
We call two links $L$ and $L'$ \emph{weakly equivalent} if $L $ is equivalent to either $L'$ or $rm(L')$.

\begin{thm}[{\cite{joy82} and \cite{Mat84}}]\label{joy}
Let $K$ and $K'$ be two knots. Then $Q(K) \cong Q(K')$ if and only if $K$ and $K'$ are weakly equivalent.
\end{thm}

\begin{defn}\label{defn:coloring}
    Let $L$ be a link and $X$ a finite quandle. A homomorphism $ Q(L) \to X$ is called a \emph{coloring} of $L$ with $X$. The number of such colorings is called the \emph{quandle counting invariant} and denoted as
    \[ Φ_X^ℤ ( L ) \coloneqq |\Qnd (Q(L), X)|.\]
\end{defn}

Again, we can think of a coloring as assigning one color, that is element of $X$ to each arc of the diagram, that is each generator of $Q(L)$. 
\dobib

\section{Fundamental quandle of linkoids}\label{fundamentalsection}
In this section we generalize the fundamental quandle to linkoids and examine how the fundamental quandle behaves under a forbidden move of endpoints of linkoids.

\begin{defn}\label{quandledeflinkoid}
    Let $L$ be an oriented linkoid diagram and $A(L)$ the set of arcs in L. The \emph{fundamental quandle} of $L$ is defined as 
    \[
       Q(L) \coloneq Q\langle x \in A(L) \quad | \quad r_τ \text{ for all crossings } τ \rangle,
    \]
    where the quandle consists of words in $A(L)$ modulo the quandle axiom relations and the relations given by each crossing as in Figure \ref{quandlecond}. 
\end{defn}

\begin{thm}
    The fundamental quandle is invariant for linkoids.
\end{thm}

\begin{proof}
    Two equivalent diagrams can be transformed into each other with a finite sequence of Reidemeister moves. Because changing a diagram by a Reidemeister move gives an isomorphic fundamental quandle, the fundamental quandle does only depend on the linkoid and not on the diagram. 
\end{proof}

Lemma \ref{linkoidcomp} was known for links, see  \cite{FR92}. We now generalize it to linkoids.

\begin{lem}\label{linkoidcomp}
    Let $L$ be a linkoid. There is a one-to-one correspondence between components of $L$ and algebraic components of $Q(L)$.
    \end{lem}
  
    \begin{proof}
       Let $a, b$ be labels of arcs in $L$. Assume $a$ and $b$ lie in the same component of $L$, then we can "walk" through the diagram from $a$ to $b$. At each crossing where we go under an arc, we need to quandle with the label of the over-crossing arc $c_i$, where $i$ counts the times we pass a crossing this way. This corresponds to the inner isomorphism $β_{c_i}$. After completing the walk from $a$ to $b$ we get a morphism $f = β_{c_n} \circ \dots \circ β_{c_1}\in Inn(Q(L))$ with $f(a) = b$.\\
       
       Now let $a$ and $b$ be two labels of arcs in $L$ and assume there is a series of generators of $Q(L)$ such that $((((a \q c_1) \q c_2) \q \dots )\q c_n) = b$. For every crossing relation $x = y \q z$ the arcs with label $x$ and $y$ lie in the same component of the linkoid and the quandle axioms do not change the element input of the inner homomorphisms, therefore $a$ and $b$ lie in the same component of $L$. \\
    
       This proves that there is a one-to-one correspondence between components of $L$ and algebraic components of $Q(L)$.
    \end{proof}
    
The main difference between links and linkoids with open components  is the existence of endpoints. We now see how the fundamental quandle acts under an under forbidden endpoints.

\begin{lem}\label{quandlefaillem}
    The fundamental quandle of a linkoid is invariant under the forbidden move $Ω_-$. 
\end{lem}

\begin{figure}[ht]
	\begin{center}
      \includegraphics[height=2cm]{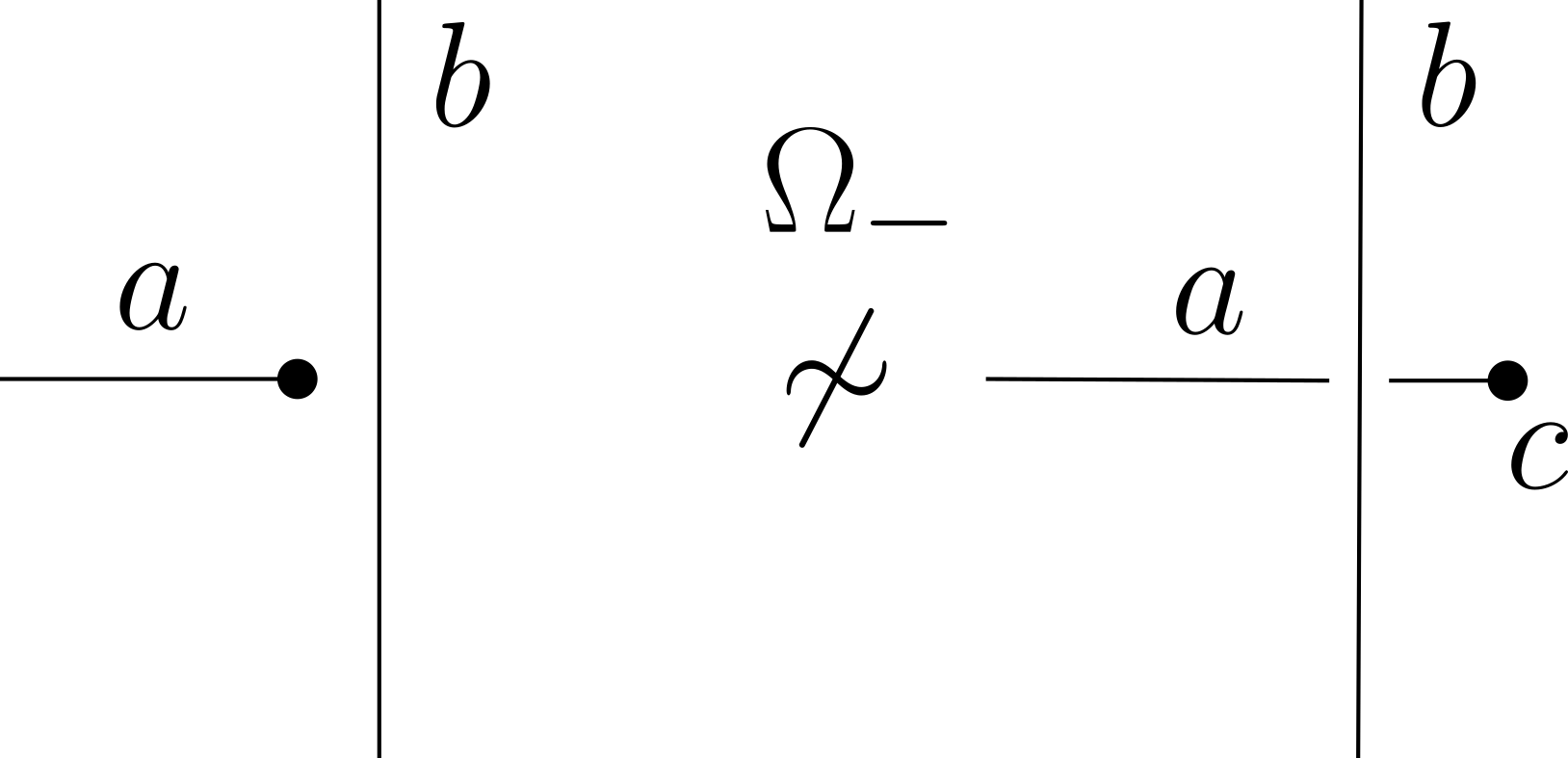}
      \caption{The under forbidden move with labeled arcs}
		\label{forbidden_l}
    \end{center}
\end{figure}

\begin{proof}
    Let $L$ and $L'$ be two linkoids that differ only by one forbidden under move $Ω_-$.  So the presentation of the fundamental quandle of one linkoid, let's say $Q(L)$, has one extra generator $c$ and one extra relation $a \q^ε b = c$ as depicted in Figure \ref{forbidden_l}. Here $ε = \pm 1$ is the sign of the crossing that  got added by the $Ω_-$-move. Note that $c$ appears in no other relation in the presentation of $Q(L)$. This means the map $φ \colon Q(L) \to Q(L')$ defined by $φ(x) = x$ for all generators $x \neq c \in Q(L)$ and $φ(c) = a \q^ε b$ (and extended to the quandle) satisfies all relations in $Q(L')$ (because the relations in both quandles are the same). Hence, $φ$ is a quandle homomorphism. Now $φ$ has the inverse map given by $φ^{-1}(x) = x \in Q(L)$, so it is an isomorphism. This shows $Q(L) \cong Q(L')$. 
\end{proof}

By Corollary \ref{omegacheat}, we deduce following.
\begin{cor}\label{quandlefail}
    If two knotoids $K$ and $K'$ represent the same knot, then $Q(K)$ is isomorphic to $Q(K')$. 
\end{cor}

This means that when we consider a coloring as a map in $\Qnd (Q(K), X)$ for some finite quandle $X$, the coloring does not depend on the specific knotoid we choose to color, but only on the knot it represents. We will enhance the fundamental quandle for linkoids in Section \ref{pointedsection} in a way that it can track the color of end-arcs.\\

On the other hand, Corollary \ref{quandlefail} means we can study knotoids to understand fundamental quandles of knot-type knotoids (and vice versa), which only differ by one relation from knot quandles. But a knotoid can potentially have fewer crossings and arcs, hence the quandle has fewer generators and relations. Specifically for computational purposes, this can be very helpful. \\

We now examine how the quandle of a knotoid corresponds to the quandle of the knot it represents. \\

First, let $K$ be a knot-type knotoid representing the knot $K_-$. Let's denote the labels of the arcs connected to the leg and head of $K$ by $l$ and $h \in Q(K)$. The presentation of the fundamental quandle of $K$ and $K_-$ only differs by the extra relation $ l = h$ in $Q(K_-)$ (because both end-arcs are connected). So we see that $Q(K_-)$ is isomorphic to $\quotient{Q(K)}{(l = h)}$, where the quotient is the quandle that is obtained by adding the new relation to the quandle $Q(K)$. 
It is then clear that if $l = h \in Q(K)$ then $Q(K) \cong Q(K_-)$. 

\begin{lem}\label{neqendpoints}
    Let $K_-$ be a knot that is not equivalent to its reversed mirror $rm(K_-)$ and let $K$ be a knot-type knotoid representing $K_-$. Denote the labels corresponding to the end arcs in $Q(K)$ by $l$ and $h$. Then $l \neq h \in Q(K)$.  
\end{lem}

\begin{rmk}
    When we talk about the mirror knotoid or mirror linkoid of a knotoid or linkoid, we have to be careful. There are two variants of mirroring a knotoid. We can either change the sign of every crossing (imagine holding the mirror behind the linkoid), or by reflecting it on an axis in the plane, outside the diagram. For links both of these are equivalent by turning the link around 180° in three-dimensional space. However, we are not allowed to do this for linkoids, due to the endpoints. Here, $m(L)$ refers to the second concept of a mirror image. By $rm(L)$ we denote this mirror linkoid with reversed orientation. 
\end{rmk}

\begin{proof}[Proof of Lemma \ref{neqendpoints}]
    By Theorem \ref{joy} we know $Q(K) \cong Q(rm(K))$. 
    The knot $K_-$ can be written as $K \cup α$ for some arc $α$. Now we take the connected sum of $K$ and $K$ or $K$ and $rm(K)$ along this arc $α$ connecting the endpoints. The composite knots $K_- \# K_-$ and  $K_- \# rm(K_-)$ can be seen as two copies of $K$ (or of $K$ and $rm(K)$) connected on the endpoints, one time head to leg and one time head to head (and reversed orientation). See Figure \ref{knotsum} below. So by assumption $rm(K_- \# K_-) \not\sim K_- \# K_- \not\sim K_- \# rm(K_-)$.\\
    
\begin{figure}[ht]
	\begin{center}
      \includegraphics[height=3cm]{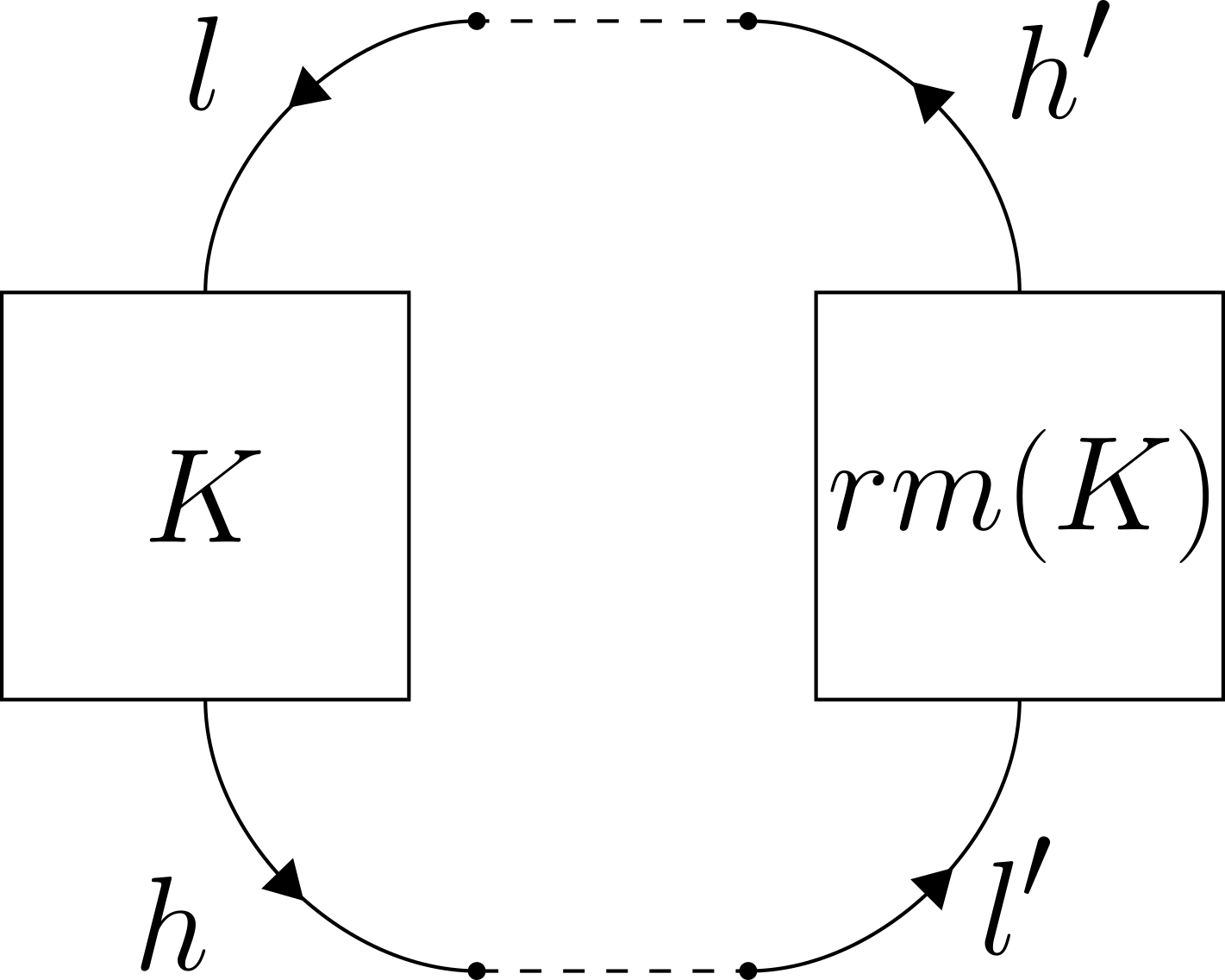}
      \caption{The knot $K_- \# rm(K_-)$ as sum of $K$ and $rm(K)$}
		\label{knotsum}
    \end{center}
\end{figure}

    We examine the fundamental quandles of the composite knots $K_- \# K_-$ and $K_- \# rm(K_-)$. Let
    \[     
        Q(K) = Q \langle l, h, x_i | r_j \rangle
    \]
    be the presentation of $Q(K)$ given by the diagram of $K$, where $l$ and $h$ are the leg and head of $K$. The generators and relations of $K_- \# K_-$ are given by
    \[
        Q(K_- \# K_-) = Q\langle l, h, x_i, l', h', x_i, x'_i | r_j, r'_j, l = l', h = h' \rangle
    \]
    where $l,h,x_i$ and $r_i$ are coming from the left and $l', h', x'_i$ and $r'_j$ right side of the composite knot. Similarly,
    \[
        Q(K_- \# rm(K_-)) = Q\langle l, h, x_i, l', h', x_i, x'_i | r_j, r'_j, l = h', h = l' \rangle.
    \]
    Now assume $l = h \in Q(K)$. Then $Q(K \# K) \cong Q(K \# rm(K))$ which contradicts the assumption that $K$ and $rm(K)$ are not equivalent by using Theorem \ref{joy}. 
 \end{proof}

    In particular $Q(K) \not\cong Q(K_-)$ in the situation of the lemma above. Together with Corollary \ref{quandlefail} we see

 \begin{cor}
    Let $K_-$ be a knot that is not equivalent to its reversed mirror $rm(K_-)$ and let $K$ be a knotoid that represents $K_-$. Then $Q(K) \not\cong Q(K_-)$.
 \end{cor}

\begin{note} Let $K$ be a knotoid that represents the trivial knot. Then we see immediately that $Q(K) \cong T_1 \cong Q(K_-)$. It is an interesting task to make the full characterization on $K$ (or its closure $K_{-}$) for $$Q(K) \cong Q(K_-).$$

\end{note}


\begin{lem}
    Let $L$ be a link-type 1-linkoid. Denote the labels corresponding to the end-arcs in $Q(L)$ by $l$ and $h$. Then $β_l = β_h$.  
\end{lem}

\begin{figure}[ht]
	\begin{center}
      \includegraphics[height=3cm]{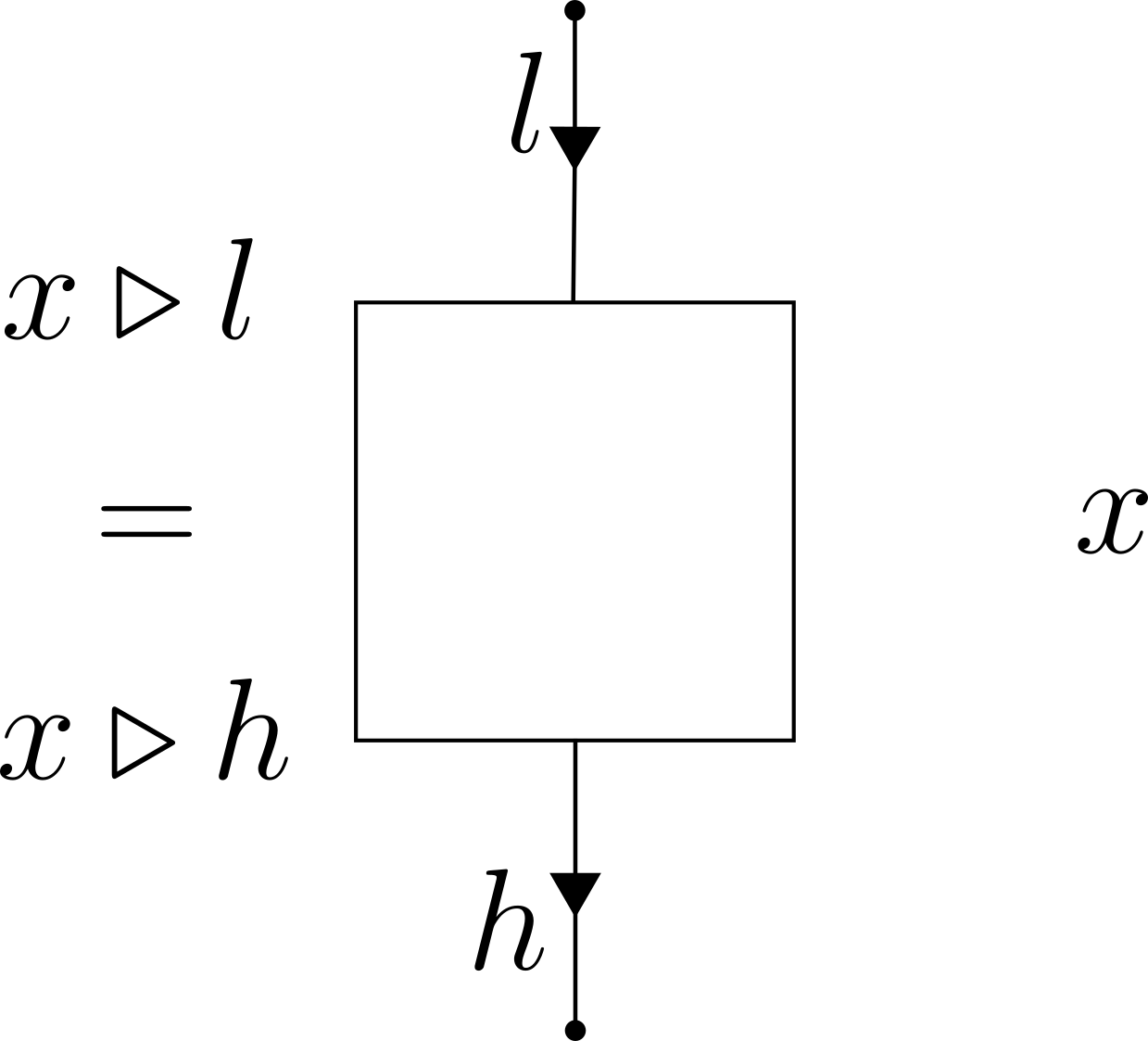}
      \caption{A link-type 1-tangle endowed with a shadow quandle coloring}
	   \label{shadowlemma}
    \end{center}
\end{figure}

\begin{proof}
    The proof can be derived by extending the argument of \cite[Lemma 5.6]{Nos11} from 1-tangles to link-type 1-linkoids. Since $L$ is a link-type $1$-linkoid, we can choose a diagram of $L$ whose endpoints lie at the same region. We consider this region to be the 'exterior' region and so the diagram of $L$ can be depicted as in Figure \ref{shadowlemma} where any closed components of the diagram lie in the box. The proof is simply based on the fact that the shadow quandle coloring is well defined.
\end{proof}

\begin{cor}
    Let $L$ be a link-type 1-linkoid. The fundamental quandle $Q(L)$ is not faithful.
\end{cor}

\begin{cor}\label{ffendmonochromatic}
    Let $L$ be a link-type 1-linkoid, $X$ a faithful quandle and $f \in \Qnd (Q(L), X)$ any coloring. Then $f(l) = f(h)$ where $l$ and $h \in Q(L)$ denote the end-arc labels. 
\end{cor}

\begin{rmk}
    Let $L$ be a link-type 1-linkoid. If for a given quandle $X$, every coloring assigns the same color to both end-arcs, the pair $(L, X)$ is called \emph{end monochromatic}. By the corollary above this is always the case if $X$ is faithful. But there are also end monochromatic pairs with non-faithful quandles. This is studied for example in \cite{CSV16} and \cite{CDS16} to better understand quandle colorings of composite knots. 
\end{rmk}

\dobib

\section{Pointed quandles}\label{pointedsection}

In this section we introduce pointed quandles and fundamental pointed quandles of ordered n-linkoids. Pointed quandles are quandles with marked points which we call basepoints. The basepoints in a fundamental pointed quandle of a linkoid encode the labels of  endpoints of the linkoid, and enables us to detect the under forbidden move on linkoids.

\subsection{n-pointed quandles}\label{intropointedsection}

\begin{defn}\label{npquandle}
    An \emph{$n$-pointed quandle} $(X, x_1, … , x_{n})$ is an ordered tuple consisting of a quandle $X$ together with $n$ (ordered) elements $x_1, … , x_{n}$ of $X$. We call $x_1, … , x_{n}$ the \emph{basepoints} of the pointed quandle. 
\end{defn}

Note that a 0-pointed quandle is again a quandle. 

\begin{defn}
    A homomorphism between two $n$-pointed quandles 	
	\[φ \colon (X, x_1, …, x_{n}) \to (Y, y_1, …, y_{n})\] 
	is a quandle homomorphism $φ \colon X \to Y$ such that $φ(x_i) = (y_i)$ for $i = 1, …, n$. We denote the category of $n$-pointed quandles by $\catname{PQnd_n}$.
\end{defn} 

There is a forgetful functor $U_n \colon \PQnd{n} \to \Qnd$ mapping $(X, x_1, …, x_n)$ to $X$, forgetting the basepoints. The set of all $n$-pointed quandles with underlying quandle $X$ is denoted by $U_n^{-1}(X)$.

For readability, we denote pointed quandles by calligraphic letters, for example $\X = (X, x_1, …, x_n)$.

\begin{rmk}\label{autoremark}
    For a given pointed quandle $\X = (X, x_1, …, x_n)$ every (unpointed) quandle homomorphism $(φ \colon X \to Y )\in \Qnd(X, Y)$ gives a pointed quandle homomorphism $\X \to (Y, φ(x_1), … , φ(x_n))$. We write 
	\[φ(\X) = (Y, φ(x_1), … , φ(x_n)).\]
	The map $φ \colon X \to Y $ is a quandle isomorphism if and only if $φ \colon \X \to φ(\X)$ is a pointed quandle isomorphism.
\end{rmk}

\begin{defn}
    An \emph{ordered $n$-linkoid diagram} is an $n$-linkoid diagram with a given ordering of the open components.
\end{defn}

Note that the ordering is invariant under Reidemeister moves. 

\begin{defn}
    Given an ordered $n$-linkoid diagram $L$, we define the \emph{fundamental pointed quandle of $L$} as the $2n$-pointed quandle 
    \[
    P(L) \coloneqq (Q(L), l_1, h_1, …, l_n, h_n)
    ,\]
    where $Q(L)$ is the fundamental quandle of $L$ as defined in \ref{quandledeflink} and $l_i$ and $h_i$ are the generators corresponding to the arcs adjacent to the leg and head of the $i$-th component of $L$. 
\end{defn}

\begin{thm}
    The fundamental pointed quandle is invariant under Reidemeister moves and the spherical move.
\end{thm}

\begin{proof}
	This follows directly from the fact that Reidemeister moves happen away from the endpoints. So the quandle isomorphism between the fundamental quandle before and after a Reidemeister move maps the endpoint labels to endpoint labels, hence it is a pointed quandle isomorphism.
\end{proof}

If $L$ is a link (that is a 0-linkoid) then the fundamental pointed quandle is a 0-pointed quandle which is simply a quandle. So $P(L) = Q(L)$. In this sense the fundamental pointed quandle is a generalization of the fundamental quandle. \\

In Example \ ref{coolpointed} we calculate the fundamental pointed quandles of two knotoids that represent the same knot via the under- closure.

\begin{figure}[ht]
	\begin{center}
        \begin{subfigure}[b]{0.4\textwidth}
            \centering
            \includegraphics[height=0.11\textheight]{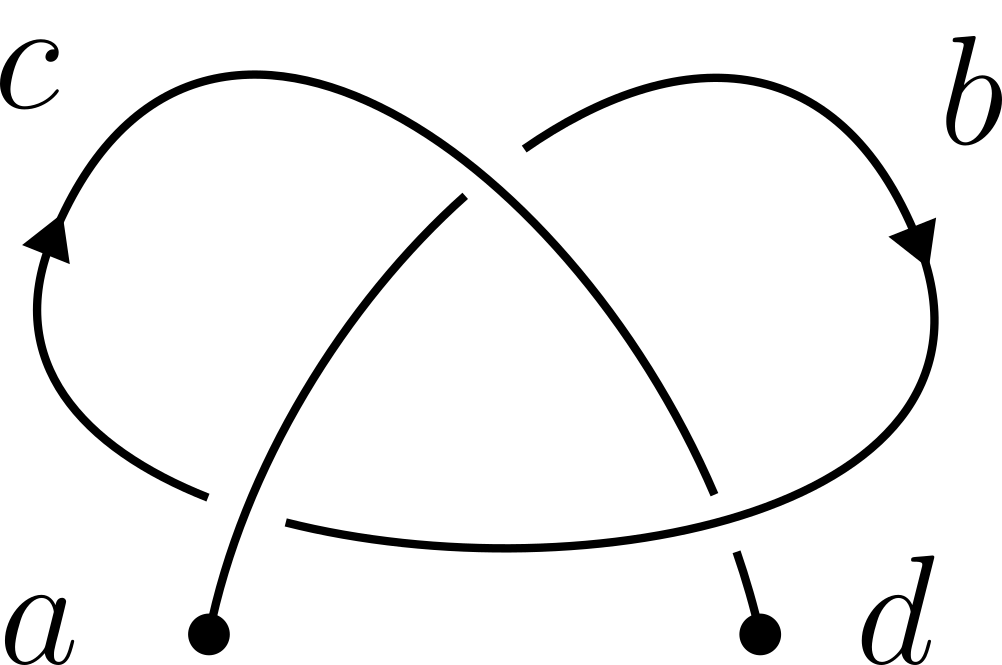}
            \caption{$K_1$}
            \label{trefoiloids1}
        \end{subfigure}
        \begin{subfigure}[b]{0.4\textwidth}
            \centering
            \includegraphics[height=0.11\textheight]{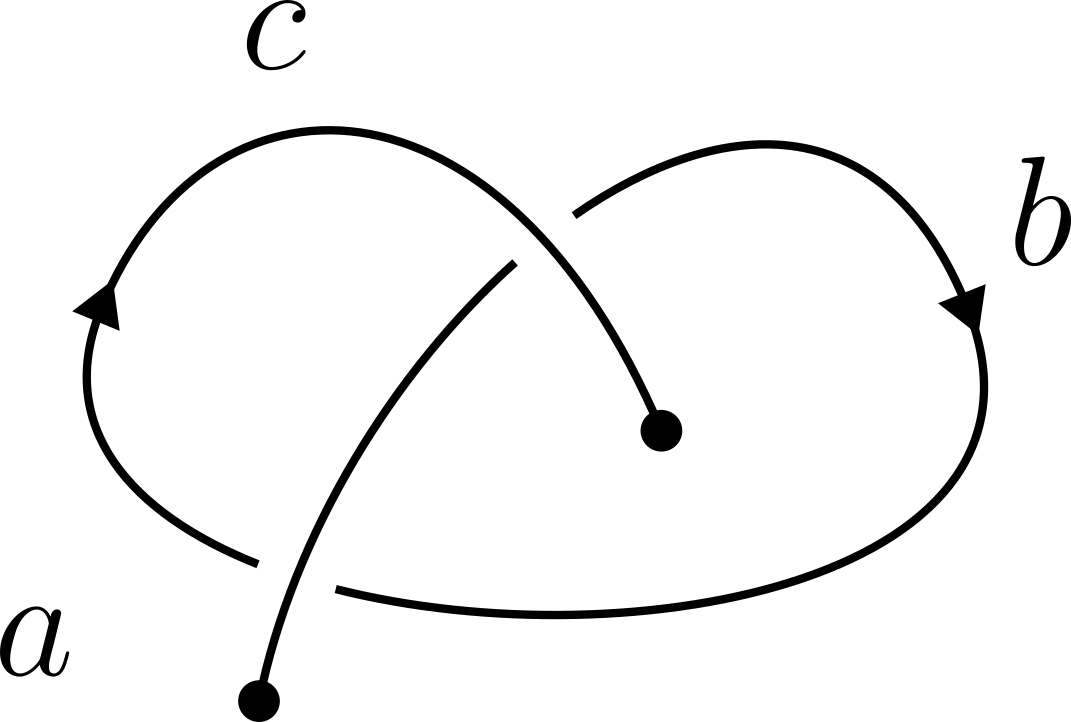}
            \caption{$K_2$}
            \label{trefoiloids2}
        \end{subfigure}
        \caption{Two knotoids representing the same knot}
        \label{trefoiloids}
	\end{center}
\end{figure}

\begin{exmp}\label{coolpointed}
    Consider the knotoids $K_1$ and $K_2$ in Figure \ref{trefoiloids}. Since $K_1$ and $K_2$ represent the trefoil knot,  their fundamental quandles, $Q(K_1)$ and $Q(K_2)$  are isomorphic by Proposition \ref{quandlefail}. We explicitly compute the quandles as
    \begin{align*}
        Q(K_1) &\cong Q\langle a, b, c, d| b = a \q c, c = b \q a , d = c \q b \rangle \\
        &\cong Q\langle a, b, c| b = a \q c, c = b \q a \rangle \cong Q(K_2).        
    \end{align*}
    Their fundamental $2$-pointed quandles are then $P(K_1) = (Q(K_1), a, c \q b)$ and $P(K_2) = (Q(K_1), a, c)$. We observe that $P(K_1)$ and $P(K_2)$  are not isomorphic pointed quandles. One can see this as follows. Let $\X = ((R_3), 0, 0)$ be a pointed dihedral quandle and consider a pointed quandle homomorphism $f \colon P(K_1) \to \X$. Because $f$ maps basepoints to basepoints we know $f(a) = 0 = f(c \q b)$. Assume $f(b)= 1$, then $f(c) = f(b \q a) = f(b) \q f(a) = 1 \q 2$. We need to check the other relation in $Q(K_1)$: $f(b) = 1 = 0 \q 2 = f(a \q 2)$. So this is a non-constant 2-pointed quandle homomorphism.\\

    Now consider a pointed quandle homomorphism $g \colon P(K_2) \to \X$. So $g(a) = 0 = f(c)$. Then $g(b) = g(a \q c) = 0 \q 0 = 0$. Hence, the constant homomorphism $g(x)= 0$ is the only such morphism. It follows that $P(K_1) \not\cong P(K_2)$.
\end{exmp}

Example \ref{coolpointed} verifies that the fundamental pointed quandle is not invariant under the under forbidden move. With the fact that it contains the information of the fundamental quandle,  the example proves that the fundamental pointed is at least powerful as the fundamental quandle for 1-linkoids that represent the same knot or link. \\

\begin{rmk}
	Note that if a knotoid $K$ represents the unknot, then $Q(K) = T_1$ has only one element. This means its fundamental pointed quandle can only have this element as basepoints. Hence, the fundamental pointed quandle of two knotoids that represent the unknot are always isomorphic. So the fundamental pointed quandle carries the same information as the fundamental quandle in this case.
\end{rmk}

Using the fundamental pointed quandle, we now set a sufficient condition for two 1-linkoids to represent the same link.

\begin{lem}\label{ltllinkequiv}
 Let $L$, $L'$ be two link-type 1-linkoids such that $P(L) \cong P(L')$ and let $L_-$ and $L'_-$ be the closure of $L$ and $L'$, respectively. Then $Q(L_-) \cong Q(L'_-)$.
\end{lem}

\begin{proof}
    Let $φ \colon P(L) = (Q(L), l, h) \to P(L') = (Q(L'), l', h')$ be the pointed quandle isomorphism. In particular, $φ$ is a quandle isomorphism. Because $φ(h) = h'$ and $φ(l) = l'$, also the map on the quotients
    \[ \tilde{φ} \colon \quotient{Q(L)}{(h = l)} \to \quotient{Q(L')}{(h' = l')}
        \]
    is a quandle isomorphism. Now  $Q(L_-) \cong \quotient{Q(L)}{(h = t)}$, which completes the proof.
\end{proof}

\begin{defn}
Let $L$ be a $1$-linkoid and $L_{\sim}$ be the link-type 1-linkoid which is derived from $L$ by  moving the head of $L$ to the region with the leg by using a finite sequence of forbidden under moves and denote its end-arcs by $\tilde{h}$ and $l \in Q(L) \cong Q(L_\sim)$. Let $f \in Inn(Q(L))$ be the map that is derived from the sequence of consecutive under forbidden moves such that $f(h)=\tilde{h}$. The relation $\tilde{h} = l$ in this setting is called the \emph{closing relation} of $L$ and denoted by $r_L$. Note that for link-type 1-linkoids the closing relation is simply $l = h$.

\end{defn}

\begin{cor} Let $L$ be a 1-linkoid, $L_{\sim}$ be the corresponding link-type 1-linkoid, and $L_{-}$ be the under-closure of $L$. Then, we have
\[
    Q(L_-) \cong \quotient{Q(L_\sim)}{(\tilde{h} = l)} \cong \quotient{Q(L)}{(f(h) = l)} = \quotient{Q(L)}{r_L}.
\]
    
\end{cor}



\begin{exmp}
Consider again the knotoid $K_2$ in Figure \ref{trefoiloids2}. Its closing relation is $r_{K_2} \colon a = c \q b = β_b (c)$ in $Q(K_2) = Q\langle a, b, c| b = a \q c, c = b \q a \rangle$, because we need to move the head under the arc labeled $b$ to move it into the region with the leg. \end{exmp}

With the closing relation we can now write the fundamental pointed quandle of $L_\sim$ as $(Q(L), l, f(h))$. Together with Lemma \ref{ltllinkequiv} this shows the following corollary.

\begin{cor}
    Let $L$ and $L'$ be two 1-linkoids with closing relations $r_L \colon f(h) = l$ and $r_{L'} \colon g(h') = l'$ as above. If there exists a pointed quandle isomorphism $φ \colon (Q(L), l, f(h)) \to (Q(L'), l', g(h'))$, then $Q(L_-) \cong Q(L'_-)$.
\end{cor}

In the case that the 1-linkoid is a knotoid, the following corollary is obtained directly by Theorem \ref{joy}.

\begin{cor}
    Let $K$ and $K'$ be two knotoids with closing relations $r_K \colon f(h) = l$ and $r_{K'} \colon g(h') = l'$. If there exists a pointed quandle isomorphism \[
		φ \colon (Q(K), l, f(h)) \to (Q(K'), l', g(h')),
	\]
	then $K$ and $K'$ represent weakly equivalent knots. 
\end{cor}

\subsection{Isomorphism classes of pointed quandles}\label{orbitsection}

For a quandle $X$, $n \in ℕ$ and the forgetful functor $U_n$, we define
\[
        P_n(X) \coloneq \quotient{U_n^{-1}(X)}{Aut(X)},
\]
the classes of pointed $n$-quandles with underlying quandle $X$ under isomorphy. Here $ \X \sim \Y \in P_n(X)$ if they are isomorphic $n$-pointed quandles. By Remark \ref{autoremark} this is exactly the case if there exists $φ \in Aut(X)$ with $\Y = φ(\X)$. We denote the number of such classes by $d_n(X) \coloneq  |P_n(X)|$.\\

Because $U_n^{-1}(X) = \{(X, x_1, …, x_n) | x_i \in X \} \cong X^n$ as sets we only write $(x_1, …, x_n)$ instead of $(X, x_1, …, x_n)$, if $X$ is clear from the context. \\

Let us compute the isomorphism classes of quandles of some examples:
\begin{exmp}\label{3qhexmp}~
    \begin{enumerate}
        \item Let $X = R_3$ be the dihedral quandle $\{0, 1, 2\}$ with $x \q y ≡ 2y - x \text{ (mod 3)}$. \\
        Note that the automorphism group $Aut(R_3)$ is a subgroup of  $S_3$, the symmetric group. This means we can describe inner automorphisms as permutations. In particular $(01) = \be{2}$, $(12) = \be{0}$ and $(02) = \be{1}$. So $S_3 = Inn(R_3) \subseteq Aut(R_3)$. \\
        In $P_1(R_3)$ we see 
        \[(R_3,0) = (0) \overset{(01)}{\sim} (1) \overset{(12)}{\sim} (2), \]
        so $P_1(R_3) = \{ [(0)] \}$ and $d_1(R_3) = 1$.\\

        For $P_2(R_3)$, we see
        \[(0, 0) \overset{(01)}{\sim} (1,1) \overset{(12)}{\sim}(2,2) \]
        and
        \[(0,1) \ \overset{(12)}{\sim} (0,2) \overset{(01)}{\sim} (1,2) \overset{(02)}{\sim} (1,0) \overset{(12)}{\sim} (2,0) \overset{(01)}{\sim} (2,1)\]
        so $P_2(R_3) = \{ [(0,0)], [(0,1)]\}$ has two equivalent classes.\\

        \item Let $X = T_3$ be the trivial quandle. We see that $Inn(T_3) = \{ id \}$, since $β_x = id$ for all $x \in T_3$. For any $f \in S_3$, we notice that $f(a \q b) = f(a) = f(a) \q f(b)$ so $ f \in Aut(T_3)$. This shows that $Aut(T_3) = S_3$ and hence $P_2(T_3) = \{ [(0,0)], [(0,1)]\}$.
        \item There are only three non-isomorphic quandles with three elements. So now let $V_3$ be the remaining quandle with three elements, that is $\be{0}= (12)$ and $\be{1} = \be{2} = id$. Let $f \in Aut(V_3)$. If $f(0) = 0$ then either $f = id$ or $f = (12)$. In both cases $f \in Inn(V_3)$. Now assume $f(0) = 1$. Then $f(2)= f(1) \q f(0) = f(1) \q 1 = f(1)$ which is a contradiction to $f$ being a bijection. Similarly, if we assume $f(0) = 2$. This shows $Aut(V_3) = Inn(V_3) = \{ id, (12) \}$. \\
        So $(1) \sim (2)$ and $P_1(V_3) = \{ [(0), (1)] \}$ with $d_1(V_3) = 2$.\\

        For $P_2(V_3)$, we observe         
        \begin{align*}
            (1,1) &\sim (2,2) \\
            (0,1) &\sim (0,2) \\
            (1,0) &\sim (2,0) \\
            (1,2) &\sim (2,1), \\
        \end{align*}
        but the other combinations are not equivalent. This gives 
        \[
            P_2(V_3) = \{[(0,0)], [(1,0)], [(0,1)], [(1,1)], [(1,2)] \}.
        \]
        In particular $d_2(V_3) = |P_2(V_3)| = 5$.
    \end{enumerate}    
\end{exmp}

Now, we study which values $d_n(X)$ can have. For a finite quandle $X$ with $k = |X|$, we immediately see $1 \leq d_n(X) \leq k^n$ because $ \{(X, x_1, …, x_n) | x_i \in X \} \cong X^n$. Our goal is to find a better lower bound for $d_n(X)$ for an arbitrary finite quandle $X$ with $k$ elements. \\

Because $Aut(X)$ is a subgroup of the symmetric group $S_k$ we will count the number of orbits of $X^n$ under the diagonal action of $S_k$. That is $(x_1, …, x_n) \sim (y_1, …, y_n)$ if there is a permutation $π \in S_k$ such that $x_i = π(y_i)$ for all $i = 1, …, n$. This is a lower bound for $d_n(X)$. We denote this number by $d_{0, n,k}$ and compute it in the following part. \\

For example $d_{0, 1,k} = 1$ for all $k$, since for every $x,y \in X$ the transposition $(xy) \in S_k$ maps $x$ to $y$, so all elements lie in the same orbit. If $k = 1$ then $d_{0, n,1} = 1$ because there is only one element in $X^n$. \\

If $k \geq 2$ then $d_{0, 2,k} = 2$ because $(x, x) \sim (y,y)$ for all $x,y \in X$ as in the case $n = 1$. To see this, we note $(x_1, x_2) \sim (y_1, y_2)$ for all $x_1 \neq x_2$ and $y_1 \neq y_2$ using the permutation $(x_2 y_2)(x_1 y_1)$ if $x_2 \neq y_1$ and the permutation $(x_1 y_1 x_2)$ if $x_2 = y_1$. But $(x, x) \not\sim (y_1, y_2)$ with $y_1 \neq y_2$ because permutations are bijections. \\

This is the case for any $n$:

\begin{lem}\label{skorbits}
    Two $n$-tuples, $x = (x_1, … , x_n)$ and $y = (y_1, …, y_n) \in X^n$, lie in the same orbit under the action of $S_k$ if and only if $x_i = x_j ⇔ y_i = y_j$ for all $i,j = 1, …, n$, i.e. they have equal entries in the same positions.
\end{lem}

\begin{proof}
    Assume $x \sim y$ so there exists $π \in S_k$ with $π(x_i) = π(y_i)$ for all $i$. Because $π$ is a bijection, we see $x_i = x_j \iff y_i = π(x_i) = π(x_j) = y_j$. \\
    On the other hand if $x_i = x_j \iff y_i = y_j$, we can define a permutation on the set $\{ x_1, … x_n, y_1, … ,y_n \}$ of unique elements in the tuples with $π(x_i) = y_i$. So $π(x) = y$ and hence $x \sim y$.
\end{proof}

Note that the statement in the lemma does not depend on $k$. By the lemma above we need to count the possibilities of how many entries are equal in a tuple and in which position these entries are. We will count these recursively over the number of entries in the tuple. \\

We extend our notation to $d_{m,n,k}$ for $m,n \in ℕ$, $k \in ℕ_{\geq 1} \cup \{ \infty \}$  and $m \leq k$ which denotes the number of equivalence classes of tuples $(x_1, …, x_{m+n}) \in X^{m+n}$ where $x_{1}, … x_{m}$ are fixed \emph{unique} elements. Of course this is not possible if we would allow $m > k$. 

For instance, $d_{1,1,k} = 2$ for $k \geq 2$ because given any fixed element $x_1$, we can either have $x_1 = x_2$ or $x_1 \neq x_2$. These are all equivalence classes. We immediately observe
\[
  d_{m, 0,k} = 1 , 
\]
since all entries are already fixed. \\

The reason we only allow unique elements in the fixed entries is that there are the same number of classes completing the tuple $(x_1, …, x_m, -, \dots, - ) \in X^{m+n}$ as there are completing the tuple $(x_1, x_1,  …, x_m, -, \dots, - ) \in X^{m+1+n}$. So while counting the orbits, we can remove a fixed element if it is already fixed in another entry.  \\ 

Now assume $k > m$. For a given a tuple $(x_1, …, x_m, -, \dots, - )$, we count the number of non-equivalent possibilities for the $(n+1)$-th entry. \\
We can either choose one of the distinct elements $x_1, …, x_m$. There are $m$ such choices. For each choice there are now $d_{m, n-1, k}$ possibilities to complete the tuple, having no new fixed element in the tuple. \\
Or we choose a new element. Then there are $d_{m+1, n-1, k}$ many ways to complete the tuple. This leads us to
\[
    d_{m, n, k} = m \cdot d_{m,n-1, k} + d_{m+1, n-1, k}.
\]
If $k = m$, then  $\{x_1, … x_m\} = X$. This means we can only choose elements that are already in the tuple. This proves the following theorem.

\begin{thm}\label{symcounti}
    Let $X$ be a set with $|X| = k \in ℕ \cup \{ \infty \}$. Let $d_{m,0,k} = 1$ and 
    \[  
    d_{m,n,k} = 
    \begin{cases}
        d_{m, n, k} = m \cdot d_{m,n-1, k} + d_{m+1, n-1, k} &\text{if } m < k \\
        d_{m, n, k} = m \cdot d_{m,n-1, k} &\text{if m = k}.
    \end{cases}
    \]
    Then $\left|\left(\quotient{X^n}{S_{k}}\right)\right| = d_{0,n,k}$.
    
\end{thm}

\begin{exmp}
    With this result we can explicitly compute $d_{m,n,k}$ for $n+m \leq k$:
    \begin{itemize}
         \item $n=1$: We see that
            \[ d_{m, 1, k} = m d_{m,0,k} + d_{m+1,0,k} = m + 1. \]
            This shows $d_{0 , 1, k} = 1$, as we have seen before. 
        \item $n = 2$:
            \[ d_{m, 2, k} = m d_{m,1,k} + d_{m+1,1,k} = m(m+1) + m+2 = m^2 + 2m + 2\]
            So as expected $d_{0 , 2, k} = 2$. 
        \item $n = 3$:
            \[ d_{m,3,k} = m d_{m,2, k} + d_{m+1,2,k} =  … = m^3 + 3m^2 + 5m + 5\]
            and $d_{0 , 3, k} = 5$.
        \item Further $d_{0 , 4, k} = 15$, $d_{0 , 5, k} = 52$ and $d_{0 , 6, k} = 203$.
    \end{itemize}

    On the other hand for $n = 3$ and $k=2$ we compute 
    \[
        d_{0,3,2} = d_{1, 2,2} = 1 \cdot d_{1,1,2} + d_{2, 1, 2} = (d_{1,0,2} + d_{2,0,2}) + 2 \cdot d_{2,0,2} = 4.
    \]
    Thinking about triplets in $\{a,b\}^3$, there are exactly the four equivalent classes $[(a,a,a)]$, $[(a,a,b)]$, $[(a,b,a)]$ and $[(b,a,a)]$. In $\{a,b,c\}^3$ there is one more class, namely $[(a,b,c)]$.
\end{exmp}

\subsection{n-homogeneous quandles}

Recall that in Section \ref{orbitsection} we defined 
\[
    d_n(X) \coloneq |P_n(X)| = \left|\left( \quotient{\{(X, x_1, …, x_n) | x_i \in X \}}{Aut(X)} \right)\right|.
\]

We give quandles with minimal $d_n(X)$ a special name and study them in more detail. 

\begin{defn}\label{nhomo}
    Let $X$ be a quandle and $n \in ℕ$. We say $X$ is \emph{$n$-homogeneous} if $d_n(X) = d_{n,0,|X|}$. We say $X$ is \emph{uniform} or \emph{$\infty$-homogeneous} if $X$ is $n$-homogeneous for all $n \in ℕ$. 
\end{defn}

\begin{prop}
    A quandle is $1$-homogeneous if and only if it is homogeneous (as defined in Definition \ref{homogeneousdef}).
\end{prop}

\begin{proof}
    Observe that $d_1(X) = 1$ means $(X, x_1) \sim (X, x_2)$ for all $x_1, x_2 \in X$. Hence, there is a quandle automorphism $f \in Aut(X)$ with $x_2 = f(x_1)$. This is exactly the definition of being homogeneous. 
\end{proof}

\begin{rmk}\label{tpremark}
In \cite{Tam13} the concept of a \emph{two-pointed homogeneous quandle} is introduced. A quandle $X$ is two-pointed homogeneous if the action of $Inn(X)$ on $U_2^{-1}(X)$ has two orbits. There we act on $U_2^{-1}(X)$ only with automorphisms that are inner automorphisms. This implies $X$ being connected. Every two-pointed homogeneous quandle is of course 2-homogeneous as in our definition above. The opposite is not true. 
\end{rmk}

\begin{lem}\label{homobasics}
    Let $X$ be a quandle, $n \in ℕ$. The following are equivalent.
    \begin{enumerate}[label=(\arabic*)]
        \item $X$ is $n$-homogeneous.
        \item $X$ is $(n-1)$-homogeneous and if $x_i \neq x_j$ and $y_i \neq y_j$ for all $i \neq j$ then $(X, x_1, …, x_n) \cong (X, y_1, …, y_n) $. 
        \item Any two $n$-pointed quandles $\mathcal{X} = (X, x_1, …, x_n)$ and $\mathcal{Y} = (X, y_1, …, y_n)$ with underlying quandle $X$ are isomorphic if and only if $x_i = x_j ⇔ y_i = y_j$ for all $i,j = 1, …, n$.
    \end{enumerate}
\end{lem}

\begin{proof}
    We first proof (1) ⇔ (3) and use this to show (2) ⇔ (3).
    \begin{itemize}
        \item[(1) ⇒ (3)] Let $X$ be $n$-homogeneous. If $\mathcal{X} \cong \mathcal{Y}$, then the right-hand side of (3) follows immediately. \\
        Let's assume $x_i = x_j ⇔ y_i = y_j$ for all $i,j = 1, …, n$ holds. By Lemma \ref{skorbits} $(x_1, …, x_n)$ and $(y_1, …, y_n)$ lie in the same orbit under the action of $S_k$. Because $Aut(X)$ is a subgroup of $S_k$, the orbit $Aut(X) \cdot x \subseteq S_k \cdot x$ for all $x \in X$. Now $X$ is $n$-homogeneous, meaning there are the same number of orbits under the action of $Aut(X)$ and $S_k$. Since orbits are either equal or disjoint the orbits must be all equal. So $Aut(X) \cdot x  = S_k \cdot x$. Hence, indeed $\X \sim_{Aut(X)} \Y \in P_n(X)$ which shows $\mathcal{X} \cong \mathcal{Y}$.
        \item[(3) ⇒ (1)] Assume (3). So $\mathcal{X} \cong \mathcal{Y}$ if and only if $(x_1,…, x_n) \sim_{S_k} (y_1, …, y_n)$. This means $S_k \cdot x = Aut(X) \cdot x$ for all $x \in X^n$. Hence, $X$ is $n$-homogeneous.
        \item[(2) ⇒ (3)] Let $\mathcal{X} = (X, x_1, …, x_n)$ and $\mathcal{Y} = (X, y_1, …, y_n)$ with $x_i = x_j ⇔ y_i = y_j$ for all $i,j = 1, …, n$. If all $x_i$ are different (and hence all $y_i$), then $\mathcal{X} \cong \mathcal{Y}$ by the second part of (2). If there are $i \neq j$ with $x_i = x_j$ (and therefore $y_i = y_j$), then $(X, x_1,…,\hat{x_j}, …, x_n) \cong (X, y_1, …, \hat{y_j}, …, x_n)$ are isomorphic $(n-1)$-pointed quandles because $X$ is $(n-1)$-homogeneous and (1) ⇔ (3). The same quandle isomorphism gives $\mathcal{X} \cong \mathcal{Y}$.
        \item[(3) ⇒ (2)] Assume (3), so $(X, x_1,… x_{n-1}, x_{n-1}) \cong (X, y_1,…, y_{n-1}, y_{n-1})$ if and only if $x_i = x_j ⇔ y_i = y_j$ for all $i,j = 1, …, n-1$. This shows that $X$ is indeed $(n-1)$-homogeneous. The second part of (2) also follows immediately. 
    \end{itemize}
\end{proof} 

In particular \ref{homobasics}(2) implies the following corollary.
\begin{cor}
	If a quandle $X$ is $n$-homogeneous, then it is $m$-homogeneous for all $m \leq n$.
\end{cor}
Lemma \ref{homobasics} above says that a quandle $X$ is $n$-homogeneous if and only if all $n$-pointed quandles $\X$ and $\Y$ with $U_n(\X) = X = U_n(\Y)$ are isomorphic whenever possible. Here possible means $x_i = x_j ⇔ y_i = y_j$ for all $i,j = 1, …, n$.

\begin{prop}\label{uniformbasics}
	Let $X$ be a quandle and $k =  |X| \in ℕ \cup \{ \infty \}$. 
	\begin{enumerate}[label=(\arabic*)]
		\item $X$ is uniform if and only if $X$ is $k$-homogeneous. 
		\item $X$ is uniform if and only if $X$ is $(k-1)$-homogeneous. 
		\item If $Aut(X) \cong S_{k}$, then $X$ is uniform.
		\item If $k < \infty$ is finite and $X$ uniform, then $Aut(X) \cong S_k$. 
	\end{enumerate}
\end{prop}

\begin{proof}
	\begin{enumerate}[label=(\arabic*)]
		\item A uniform quandle is obviously $k$-homogeneous. Let $X$ be $k$-homogeneous. If $k = \infty$, we are done. So let $k \in ℕ$ be finite. By Lemma \ref{homobasics} (2), it is $n$-homogeneous for all $n \leq k$. On the other hand, if $X$ is $n$-homogeneous with $n \geq k$, there is no pointed $(n+1)$ quandle $(X, x_1, … x_{n+1})$ with  $y_i \neq y_j$ for all $i \neq j \in \{1, …, n+1 \}$. So again Lemma \ref{homobasics} (2) shows that $X$ is $(n+1)$-homogeneous. Because $X$ is $k$-homogeneous, it is by induction uniform.
		\item Let $X$ be a $(k-1)$-homogeneous quandle. We want to show that $X$ is $k$-homogeneous and use (1). Again by Lemma \ref{homobasics} it is enough to show $(x_1, …, x_k) \sim (y_1, …, y_k)$ for $x_i \neq x_j$ and $y_i \neq y_j$ for all $i \neq j$. We know that there exists $f \in Aut(X)$ with $f(x_i) = y_i$ for $i = 1, …, k-1$ since $X$ is $(k-1)$-homogeneous. Because $k = |X|$ and $f$ is a bijection we find $f(x_k) = y_k$, so indeed $(x_1, …, x_k) \overset{f}{\sim} (y_1, …, y_k)$.
		\item Follows immediately from the definition together with \ref{symcounti}.
		\item Let $X= \{x_1, …, x_k\}$ be uniform and assume $Aut(X) \subsetneqq S_k$. Then there exists $f \in S_k \setminus Aut(X)$. This means $(x_1, …, x_k) \underset{Aut(X)}{\nsim} (f(x_1), …, f(x_k))$ by assumption. This contradicts $X$ being $k$-homogeneous because of \ref{homobasics}(2).
	\end{enumerate}
\end{proof}
 
\begin{exmp}
	\begin{itemize}~
		\item All trivial quandles are uniform, since $Aut(T_k) = S_{k}$ for all $k \in ℕ \cup \{ \infty \}$. Note that trivial quandles are not two-point homogeneous as described in Remark \ref{tpremark}. 
		\item $R_3$  is uniform as we have seen in Example \ref{3qhexmp}.
		\item $V_3$ in Example \ref{3qhexmp} is not 2-homogeneous. 
		\item Let us consider the regular tetrahedron quandle $X = \{0, 1, 2, 3 \}$ with 
		\[
			 β_0 = (123), \qquad β_1 = (032), \qquad β_2 = (013), \qquad β_3 = (021).
		\]
		Note that $X$ is connected, so in particular 1-homogeneous. Now $(x,y) \overset{β_x^i}{\sim} (x,z)$ for some $i$ for all $y \neq x \neq z$. This way we can reach every tuple which does not have equal entries. Here we see for example
		\[(0,1) \overset{β_0}{\sim} (0,2) \overset{β_2}{\sim} (1,2) \overset{β_1}{\sim} (1,0). \]
		So $X$ is 2-homogeneous. 
		Assume there exists $f \in Aut(X)$ such that $(0,1,2) \overset{f}{\sim} (1,0,2)$, so $1 = f(0) = f(2 \q 1) = f(2) \q f(1) = 2 \q 0 = 3$ which is a contradiction. Hence, $(0,1,2) \nsim (1,0,2)$. This means $X$ is not 3-homogeneous and in particular not uniform.
	\end{itemize}
\end{exmp}
	
The last example can be generalized to a wider set of examples. A finite quandle $X$ with $k = |X|$ is called of \emph{cyclic type} if for every $x \in X$, the map $β_x$ acts on $X \setminus \{x\}$ as a cyclic permutation of order $(k-1)$ as defined in \cite[Def. 3.5]{Tam13}. We immediately see that for example the regular tetrahedron quandle is of cyclic type. 

\begin{prop}[{\cite[Prop. 3.6]{Tam13}}]
	Every finite quandle of cyclic type is 2-homogeneous. 
\end{prop}

Now we classify all uniform quandles.
\begin{thm}\label{alluniforms}
	Let $X$ be a uniform and finite quandle. Then either $X \cong R_3$ or $X$ is trivial. 
\end{thm}

\begin{proof}
	Let $X$ be a nontrivial uniform quandle. This means there is $a_0, a_1, a_2 \in X$ with $ a_0 \q a_1 = a_2$ and $a_0 \neq a_2$. Then also $a_0, a_2 \neq a_1$. Since $X$ is uniform, $Aut(X) \cong S_{|X|}$ by Proposition \ref{uniformbasics}, in particular are the transpositions $(a_i a_j)$ in $Aut(X)$. So $a_1 = (a_1 a_2)(a_2) = (a_1 a_2) (a_0 \q a_1) = a_0 \q a_2$, similarly $a_0 = a_2 \q a_1$ and therefore 
	\[a_i \q a_j = a_k = (a_i a_j)(a_k) = (a_i a_j) (a_i \q a_j) = a_j \q a_i\]for all $\{ i,j,k\} = \{ 0,1,2\}$. 
	This computes the following section of the operation table of $X$.
	\[
	\begin{array}{c|c c c c }
	\q & a_0 &a_1& a_2 & … \\
	\hline
	a_0 & a_0 & a_2 & a_1 &  \\
	a_1 & a_2 & a_1 & a_0 & ? \\
	a_2 & a_1 & a_0 & a_2 & \\
	\vdots & & ? & & 
	\end{array}
	\]
	If $X$ has only three elements, this means it is isomorphic to $R_3$. \\
	Suppose $X$ has another element $x \neq a_0, a_1, a_2$ and assume $a_i \q x \neq a_i$ for some $i$. Then choose $a_j \neq a_i, a_i \q x$ (which exists since we have three elements to choose from). Now we see 
	\[ a_i \q x \underset{a_i \q x \neq a_i, a_j}{=} (a_ia_j)(a_i \q x) = a_j \q x, \] 
	which is a contradiction to $β_x$ being bijective. So $a_i \q x = a_i$ for all $i$. \\
	But now we see, using the transposition $(a_0x) \in Aut(X)$, that $(a_0x)(a_1 \q x) = (a_0x)(a_1) = a_1$ but $(a_0x)(a_1) \q (a_0x)(x) = a_1 \q a_0 = a_2$. But this is a contradiction to $(a_0x)$ being a homomorphism. This shows that $X$ either has only three elements or is trivial. 
\end{proof}

For $x \in X$, denote $Aut(X)_x \coloneqq \{f \in Aut(X)| f(x)= x \}$ the stabilizer subgroup with respect to $x$. We follow the proof of \cite[Proposition 3.3]{Tam13} to prove the next Proposition.
\begin{prop}
	Let $X$ be a quandle with $|X| \geq 3$. Then the following are equivalent:
	\begin{enumerate}[label=(\arabic*)]
		\item $X$ is 2-homogeneous.
		\item For every $x \in X$, the action of $Aut(X)_x$ on $X \setminus \{x \}$ is transitive.
		\item $X$ is homogeneous and there exists $x \in X$ such that the action of $Aut(X)_x$ on $X \setminus \{ x \}$ is transitive.
	\end{enumerate}
\end{prop}

\begin{proof}~\\
	\begin{enumerate}[label=(\arabic*)]
		\item[(1) ⇒ (2)] Let $X$ be 2-homogeneous. For arbitrary $x \in X$, we know that $(x, y_1) \overset{f}{\sim} (x, y_2)$ for any $y_1, y_2 \neq x$ for some $f \in Aut(X)$. So $f \in Aut(X)_x$ and $f(y_1) = y_2$. This shows that the action of $Aut(X)_x$ is transitive on $X \setminus \{x \}$.
		\item[(2) ⇒ (3)] We only need to show that $X$ is homogeneous. For any $y_1,y_2 \in X$ we can choose $x \in X\setminus \{y_1, y_2 \}$ because $ |X| \geq 3$. By assumption, there exists $f \in Aut(X)_x \subseteq Aut(X)$ with $f(y_1) = y_2$. 
		\item[(3) ⇒ (1)]  Assume (3). So $X$ is 1-homogenous. By Lemma \ref{homobasics}, it is enough to show that $(x_1, x_2) \sim (y_1, y_2)$ for all $x_1 \neq x_2$ and $y_1 \neq y_2$. By assumption, there exists $x \in X$ such that $Aut(X)_x$ acts transitive on $X \setminus \{x\}$. Because $X$ is (1-)homogeneous, there exists $f,g \in Aut(X)$ with $f(x_1) = x$ and $g(y_1) = x$. Because $x_1 \neq x_2$ and $y_1 \neq y_2$ and $f$ and $g$ are bijective, we see $f(x_2), g(y_2) \neq x$. Hence, there exists $h \in Aut(X)_x$ with $h(f(x_2)) = g(y_2)$. Now we can see
		\[ g^{-1} \circ h \circ f (x_1, x_2) = g^{-1} \circ h (x, f(x_2)) = g^{-1}(x, g(y_2)) = (y_1, y_2)
		\]
		which shows $(x_1, x_2) \sim (y_1, y_2)$. 
		
	\end{enumerate}
\end{proof}

\subsection{Pointed quandle counting invariant and quandle counting matrix}

We now turn our attention back towards linkoids and introduce $n$-pointed quandle colorings of linkoids.

\begin{defn}\label{pointedcoloring}
	Let $\X$ be a $2n$-pointed quandle and $L$ an $n$-linkoid. We define the \emph{pointed quandle counting invariant} as 
	\[
		Φ_{\X}^ℤ(L) \coloneqq | \PQnd{2n} (P(L), \X)|.	
	\]
\end{defn}

\begin{thm}
	The pointed quandle counting invariant is invariant under Reidemeister moves.
\end{thm}

\begin{proof}
	This follows immediately from the fact that $P(L)$ is  a linkoid invariant.
\end{proof}

Note that the (unpointed) quandle counting invariant with respect to a quandle $X$ is always at least 1 because we can color the diagram trivially. For the pointed quandle counting invariant this is in general not true. If not all basepoints of a pointed quandle $\X$ are equal, we cannot color the diagram trivially. This means $Φ_{\X}^ℤ(L)$ can be zero.\\

We combine all pointed quandles in $U_2^{-1}(X)$, the set of 2-pointed quandles with underlying quandle $X$, into one matrix.

\begin{defn}\label{pqcmatrix}
	Let $L$ be a 1-linkoid and $X = \{1, …, k\}$ a finite quandle. We define the \emph{quandle counting matrix} $Φ_X^{M_k}(L)$ of $L$ with respect to $X$ by the $k \times k$ matrix
	\[
	\left(Φ_X^{M_k}(L)\right)_{i,j} \coloneqq Φ_{(X,i,j)}^ℤ (L),
	\]
	with the pointed quandle counting invariant for each possible combination of basepoints as entries.
\end{defn}

\begin{thm}
	Let $X$ be a finite quandle, $L$ a 1-linkoid. The quandle counting matrix  $Φ_{(X,i,j)}^{M_k} (L)$ is an invariant of $L$.
\end{thm}

\begin{proof}
	This follows immediately from the fact that the pointed quandle counting invariant is invariant under Reidemeister moves. 
\end{proof}

We first collect some basic properties and observations and then see some examples. Therefore, we use the fact that for two isomorphic $2n$-pointed quandles $\X$ and $\Y$ the pointed quandle counting invariant is equal, that is $Φ_\X^ℤ(L) = Φ_\Y^ℤ(L)$ for any $n$-linkoid. \\

\begin{prop}\label{pqcmatrixprop}
	Let $L$ be a 1-linkoid and $X$ a finite quandle, $i,j \in X$. Then 
	\begin{enumerate}[label=(\arabic*)]
		\item $(Φ_X^{M_k}(L))_{i,j} \geq 0$.
		\item $(Φ_X^{M_k}(L))_{i,i} \geq 1$.
		\item $Φ_X^{M_k}(L) = I_k$ the identity matrix if and only if $L$ is only trivially colorable by $X$.
		\item $\sum_{i,j = 1}^{k} (Φ_X^{M_k}(L))_{i,j} = Φ_X^ℤ(L) = | \catname{Qnd}(Q(L), X) |$ the (unpointed) quandle coloring counting invariant. 
		\item If $X$ is homogeneous, then $(Φ_X^{M_k}(L))_{i,i} = (Φ_X^{M_k}(L))_{j,j}$.
		\item If $X$ is 2-homogeneous, then $(Φ_X^{M_k}(L))_{i_1,j_1} = (Φ_X^{M_k}(L))_{i_2,j_2}$ for all $i_1 \neq j_1$ and $i_2 \neq j_2 \in X$. 
		
	\end{enumerate}
\end{prop}

\begin{proof}
	(1) is obvious, (2) and (3) follow from the fact that every knotoid is trivially colorable. (4) is the fact that the quandle coloring counting invariant counts all possible colorings, independent of the endpoint colors. \\
	To see (5) and (6), note that if $(i_1, j_1) \sim (i_2, j_2) \in \quotient{X^2}{Aut(X)} \cong P_2(X)$, then 
	$Φ_{(X,i_1,j_1)}^ℤ (L) = Φ_{(X,i_2,j_2)}^ℤ (L)$ since the pointed quandles are isomorphic. So if $X$ is homogeneous then $(i,i) \sim (j,j)$ for all $i,j \in X$ and if $X$ is 2-homogeneous $(i_1, j_1) \sim (i_2, j_2)$ for all $i_1 \neq j_1$ and $i_2 \neq j_2 \in X$.\\ 
	
\end{proof}
Note that any 2-homogeneous quandle is also homogeneous, so all values on the diagonal are also equal, but not necessarily equal to the values not on the diagonal. \\

\begin{lem}\label{algcomponentspqcm}
	 Let $L$ be a 1-linkoid and $X$ a finite quandle. If $i,j \in X$ are in the same algebraic component of $X$ (that is in the same orbit under the action of $Inn(X)$ on $X$), then $(Φ_X^{M_k}(L))_{i,i} = (Φ_X^{M_k}(L))_{j,j}$.
\end{lem}

\begin{proof}
	If $i,j$ lie in the same component of $X$ then there exists an inner automorphism $f \in Inn(X) \subseteq Aut(X)$ with $f(i) = j$. So $(i,i) \overset{f}{\sim} (j,j)$. 
\end{proof}

\begin{prop}\label{2pqcprop}
	Let $L$ be a link-type 1-linkoid and $X$ a finite quandle, $i,j \in X$. Then: 
	\begin{enumerate}[label=(\arabic*)]
		
		\item The trace $tr\left(Φ_X^{M_k}(L) \right) = Φ_X^ℤ (L_{-})$, the (unpointed) quandle counting invariant of the under-closure of $L$. \\
		\item If $X$ is faithful, then $Φ_X^{M_k}(L))_{i,j} = 0$ for all $i \neq j$
	\end{enumerate}
\end{prop}

\begin{proof}
	To see (1), we note that  
	\[ \Qnd( Q(K_-), X ) \cong \Qnd\left( \quotient{(Q(K))}{(h = l)}, X \right) \cong \bigcup_{i = 1}^k \PQnd{2}(P(K), (X, i,i))\] as sets. The first bijection follows from the fact that $Q(K_-) \cong \quotient{Q(K)}{(h=l)}$ where $h$ and $l$ are the arcs attached to the head and leg of $X$. To see the second bijection we note that a quandle homomorphism $φ \colon \quotient{(Q(K))}{(h = l)} \to X$ is the same as a quandle homomorphism $Q(K) \to X$ with $φ(h) = φ(l)$ and therefore an element in $\PQnd{2}(P(K), (X, φ(l),φ(h)))$. And of course every pointed homomorphism in $ψ \in \PQnd{2}(P(K), (X, i,i)$ satisfies $ψ(h) = i = ψ(l)$ and so is an element of $\Qnd\left( \quotient{(Q(K))}{(h = l)}, X \right)$. \\
	(2) follows immediately from Lemma \ref{ffendmonochromatic} which implies that every coloring of a link-type 1-linkoid by a faithful quandle assigns the same color to both endpoints. Hence, $Φ_X^{M_k}(K)_{i,j} = 0$ for all $i \neq j$.
\end{proof}

Combining Proposition \ref{pqcmatrixprop} (4) and \ref{2pqcprop} proves the following corollary.

\begin{cor}
		Let $L$ be a link-type 1-linkoid and $X$ a finite, faithful quandle. Then $Φ_X^ℤ(L) = Φ_X^ℤ(L_-)$.
\end{cor}

We now compute some quandle counting matrices. 

\begin{exmp}
	Let $K_1$ be the knotoid from Figure \ref{trefoiloids}. So 	
	\begin{align*}
    Q(K_1) &\cong Q\langle a, b, c, d| b = a \q c, c = b \q a , d = c \q b \rangle \\
     &\cong Q\langle a, b, c| b = a \q c, c = b \q a \rangle     
    \end{align*} 
	and $P(K_1) = (Q(K_1), a, c \q b)$. \\
	
	Consider the pointed quandle $X = (R_3, 0,0)$. A homomorphism $P_2(K_1) \to X$ maps $a,d \mapsto 0$. For any given value of $f(b)$ we compute $f(c) = f(b \q a) = f(b) \q 0 =  2 \cdot 0 - f(b) = - f(b) \text{ (mod 3)} $. We need to check the other relations in $Q(K_1)$, to see which values give indeed a quandle homomorphism. That is $f(a \q c) = 0 \q f(c) = 2 \cdot 0 - f(c) = f(c)  \text{ (mod 3)}$ and $f(c \q b) = f(c) \q f(b) = 2 f(b) - f(c) = 3f(b) = 0  \text{ (mod 3)}$. So every value for $f(b)$ determines exactly on coloring by $\X$. This shows $Φ^ℤ_X(K_1) = 3$. \\
	
	Now we can compute the quandle coloring matrix of $K_1$. Since $R_3$ is 2-homogeneous, all diagonal entries are equal by Proposition \ref{pqcmatrixprop}. By Proposition \ref{2pqcprop} all other entries are zero, since $K_1$ is of knot type. In total, this gives
	\[ Φ_{R_3}^{M_3} (K_1) = 
	\begin{bmatrix}
	3 & 0 & 0 \\
	0 & 3 & 0 \\
	0 & 0 & 3
	\end{bmatrix}.\]
	
	Computing the trace $tr( Φ_{R_3}^{M_3} (K_1)) = 9 = Φ_{R_3}^ℤ (K_-)$, we again obtain the usual quandle counting invariant, as expected.\\
	
	Let now $K_2$ be the other knotoid in Figure \ref{trefoiloids}. Because $Q(K_1) \cong Q(K_2)$ we write $P(K_2) = (Q(K_1), a, c)$. A homomorphism $f \colon P_2(K_2) \to (R_3, i, j)$ maps $a \mapsto i$ and $c \mapsto j$. So $f(b) = f(a \q c) = i \q j$ hence $f$ is already completely determined. Because $f(b \q a) = (i \q j) \q i = 2i - (i \q j) = 2i - 2j + i ≡ j =  f(c) \text{ (mod 3)}$ every such map is a quandle homomorphism. Then $Φ_{(R_3, i, j)} ^ℤ (K_2) = 1$ for all $i,j$ and 
	\[
	 Φ_{R_3}^{M_3} (K_2) = 
	\begin{bmatrix}
	1 & 1 & 1 \\
	1 & 1 & 1 \\
	1 & 1 & 1
	\end{bmatrix}.
	\]
	
	\begin{figure}[ht]
		\begin{center}
			\includegraphics[width=0.3\linewidth]{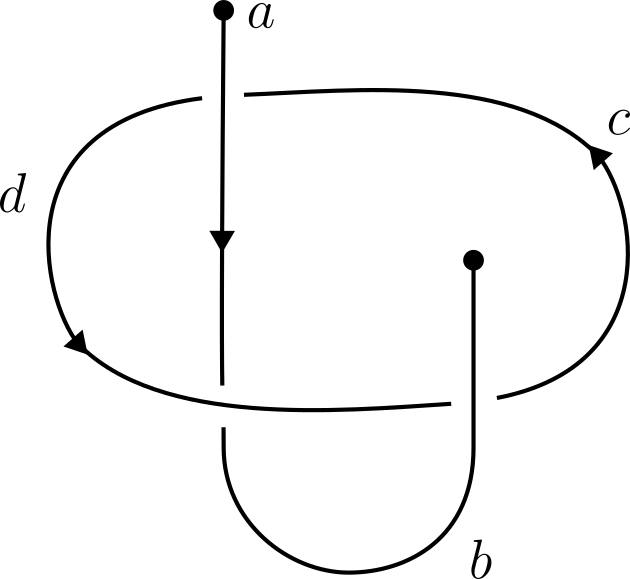}
			\caption{}
			\label{1linkoidlabeled}
		\end{center}
	\end{figure}
	As a last example consider the 1-linkoid $L$ in Figure \ref{1linkoidlabeled}
	and the quandle $V_3$ with three elements and $\be{0}= (12)$ and $\be{1} = \be{2} = id$. We have seen in Example \ref{3qhexmp} that $d_2(V_3) = 5$. This means there are five isomorphism classes of pointed quandles with underlying quandle $V_3$. So we have to compute five pointed coloring counting invariants for the coloring matrix. \\
	The fundamental quandle of $L$ is given by
	\[ Q(L) = Q\langle a, b, c, d | c = d \q b, a = b \q c, d = c \q a \rangle. \]
	We now compute the pointed quandle counting invariant for all five isomorphism classes of pointed quandles with underlying quandle $V_3$. 
	\begin{itemize}
		\item A pointed quandle homomorphism $f \colon P(L) \to (V_3, 0, 0)$ maps $a, b \mapsto 0$. If $f(c) = 0$ then $f(d) = f(0 \q 0) = 0$ gives only the trivial coloring. If $f(c) = 1$ then $f(d) = f(c \q a) = 1 \q 0 = 2$. This satisfies $f(d \q b) = 2 \q 0 = 1 = f(c)$ and $f(b \q c) = 0 \q 1 = 0 = f(a)$. So  $Φ_{(V_3, 0, 0)} ^ℤ (L) = 3$.
		\item Now consider a map $f \colon P(L) \to (V_3, 0, 1)$. Then $0 = f(a) = f(b \q c) = 1 \q f(b)$ which cannot happen. 
		\item Similarly, for a map $f \colon P(L) \to (V_3, 1, 0)$ we have $0 = f(b) = f(a \q^{-1} c) = 1 \q^{-1} f(b)$ which is also not possible.
		\item A homomorphism $f \colon P(L) \to (V_3, 1, 1)$ map $f(c) = f(d \q b) = f(d) \q b = f(d)$. Now $f(a) = f(b) \q f(c)$ gives $1 = 1 \q f(c)$ so $f(c) = 1,2$. Both values satisfy the third relation $f(d) = f(c \q a) = fa(c)$. This shows $Φ_{(V_3, 0, 0)} ^ℤ (L) = 2$.
		\item  Lastly let $f \colon P(L) \to (V_3, 1, 2)$. Because $1 = f(a) = f(b) \q f(c) = 2 \q f(c)$ we obtain $f(c) = 2$. Then $f(d) = 2 \q 1 = 2$. This satisfies $f(d \q b) = 2 \q 1 = 2 = f(c)$. So $Φ_{(V_3, 1, 2)} ^ℤ (L) = 1$. 
	\end{itemize}
	This yields the quandle counting matrix 
	\[
		Φ_{V_3}^{M_3} (L) = 
	   \begin{bmatrix}
	   3 & 0 & 0 \\
	   0 & 2 & 1 \\
	   0 & 1 & 2
	   \end{bmatrix}.
	   \]
	
\end{exmp}

\section{Discussion and further directions}\label{discussionsection}

    Recall that \emph{knot} quandles can be defined as homotopy classes of paths  from a chosen basepoint to the knot (or to a tubular neighborhood of the knot) as described in \cite[Chapter 4.4]{en15}. In this setting the expression $x \q y$ means to first walk a loop around the arc corresponding to $y$ and then walk (the path to) $x$. \\

    Knotoids can be identified with $\theta$-curves in $S^3$ \cite{Tur10}. However, the fundamental quandle of this spatial graph does not coincide with the knotoid quandle of the corresponding knotoid. It is an interesting task to find out what a knotoid quandle evaluates geometrically. \\

    In \cite{Pfl23}, the quandle 2-cocycle invariant is discussed. Unlike for knots, the 2-cocyle invariant with respect to a 2-coboundary is in general not zero for knotoids. This implies, the invariant depends on the explicit cocycle and not only its cohomology class. 
    It would be promising to present a cohomology based knotoid invariant. We expect that this can be done by utilizing pointed quandles. \\

    The relation between the fundamental quandle and the fundamental biquandle of a  link (that may be a classical or a surface link) is studied in \cite{IT20}. It is an interesting task to examine the relation between the fundamental quandle and the fundamental biquandle of a linkoid which has also introduced in \cite{Pfl23}. We expect that this relation results in a better algebraic understanding of the fundamental pointed quandles of linkoids by using the biquandle theory.

    \dobib

\printbibliography

\end{document}